\documentclass[a4paper,12pt]{article}
\begin{document}

\title{Unbounded normal operators in octonion Hilbert spaces and their spectra}
\author{Ludkovsky S.V.}
\date{15 February 2012}
\maketitle

\begin{abstract}
Affiliated and normal operators in octonion Hilbert spaces are
studied. Theorems about their properties and of related algebras are
demonstrated. Spectra of unbounded normal operators are
investigated. \footnote{key words and phrases: non-commutative
functional analysis, hypercomplex numbers, quaternion skew field,
octonion algebra, operator, operator algebra, spectra, spectral
measure, non-commutative integration
 \\
}

\end{abstract}

\section{Introduction}
Unbounded normal operators over the complex field have found
many-sided applications in functional analysis, differential and
partial differential equations and their applications in the
sciences \cite{danschw,jungexu10,kadring,killipsimon09,zelditch09}.
On the other hand, hypercomplex analysis is fast developing,
particularly in relation with problems of theoretical and
mathematical physics and of partial differential equations
\cite{brdeso,gilmurr,guesprqa}. The octonion algebra is the largest
division real algebra in which the complex field has non-central
embeddings \cite{dickson,baez,kansol}. It is intensively used
especially in recent years not only in mathematics, but also in
applications \cite{emch,guetze,girard,krausryan,kravchot}.
\par In previous works analysis over quaternion and octonions was
developed and spectral theory of bounded normal operators and
unbounded self-adjoint operators was described \cite{ludoyst,ludfov,
lujmsalop,lufjmsrf,ludancdnb}. Some results on their applications in
partial differential equations were obtained
\cite{luspraaca,ludspr,ludspr2,ludvhpde,ludcmft12}. This work
continuous previous articles and uses their results. The present
paper is devoted to unbounded normal operators and affiliated
operators in octonion Hilbert spaces, that was not yet studied
before.
\par Frequently in practical problems, for example, related with partial differential operators,
spectral theory of unbounded normal operators is necessary. This
article contains the spectral theory of unbounded affiliated and
normal operators. Notations and definitions of papers
\cite{ludoyst,ludfov,lujmsalop,lufjmsrf,ludancdnb} are used below.
The main results of this paper are obtained for the first time.

\section{Affiliated and normal operators}
\par {\bf 1. Definitions.} Let $X$ be a Hilbert space over the
Cayley-Dickson algebra ${\cal A}_v$, $~2\le v$, and let $\sf A$ be a
von Neumann algebra contained in $L_q(X)$. We say that an operator
$Q\in L_q(X)$ quasi-commutes with $\sf A$ if the algebra $alg_{{\cal
A}_v} (Q,B)$ over ${\cal A}_v$ generated by $Q$ and $B$ is
quasi-commutative for each $B\in \sf A$. A closed $\bf R$
homogeneous ${\cal A}_v$ additive operator $T$ with a dense ${\cal
A}_v$ vector domain ${\cal D}(T)\subset X$ is said to be affiliated
with $\sf A$, when \par $(1)$ $U^*TUx=Tx$ \\
for every $x\in {\cal D}(T)$ and each unitary operator $U\in L_q(X)$
quasi-commuting with $\sf A$.  The fact that $T$ is affiliated with
$\sf A$ is denoted by $T\eta \sf A$. \par An ${\cal A}_v$ vector
subspace $\bigcup_n \mbox{ }_nF{\cal D}(T)$ is called a core of an
operator $T$, if $\mbox{}_nF$ is an increasing sequence of ${\cal
A}_v$ graded projections and $\bigcup_n \mbox{ }_nF{\cal D}(T)$ is
dense in ${\cal D}(T)$.

\par {\bf 2. Note.} Definition 1 implies that $U{\cal D}(T)={\cal
D }(T)$. If $V$ is a dense ${\cal A}_v$ vector subspace in ${\cal
D}(T)$ and $T|_V\eta \sf A$, then $T\eta \sf A$. Indeed, $\lim_n
Uy^n=Uy$ for each unitary operator $U\in L_q(X)$ and every sequence
$y^n \in V$ converging to a vector $y\in {\cal D}(T)$ and with
$\lim_n Ty^n=Ty$. Therefore, the limit $\lim_n TUy^n=  \lim_n U T
y^n= UTy$ exists. The operator $T$ is closed, hence $Uy\in {\cal
D}(T)$ and $UTy= TUy$. Thus ${\cal D}(T)\subset U^* {\cal D}(T)$.
The proof above for $U^*$ instead of $U$ gives the inclusion ${\cal
D}(T)\subset U {\cal D}(T)$, consequently, $U{\cal D}(T)={\cal
D}(T)$ and hence ${\cal D}(U^*TU)={\cal D}(T)$ and $TUy=UTy$ for
every $y\in {\cal D}(T)$.
\par For an $\bf R$ homogeneous ${\cal A}_v$ additive (i.e.
quasi-linear) operator $A$ in $X$ with an ${\cal A}_v$ vector domain
${\cal D}(A)$ the notation can be used:
\par $(1)$ $Ax=\sum_j A^{i_j}x_j$ for each \par $(2)$ $x=\sum_j x_ji_j\in {\cal D}(A)\subset X$
\quad with $x_j\in X_j$ and \par $(3)$ $A^{i_j}x_j := A(x_ji_j)$ for
each $j=0,1,2,...$. \par That is $A^{i_j}(i_j^*\pi ^j) = A\pi ^j$,
where $\pi ^j: X\to X_ji_j$ is an $\bf R$ linear projection with
$\pi ^j(x)=x_ji_j$ so that $\sum_j \pi ^j =I$.
\par It can be lightly seen, that Definition 1 is natural. Indeed, if
$\sf A$ is a quasi-commutative
algebra over the Cayley-Dickson algebra ${\cal A}_v$ with $2\le v$,
then an algebra ${\sf A}_0i_0\oplus {\sf A}_ki_k$ is commutative for
$k\ge 1$ over the complex field ${\bf C}_{i_k} := {\bf R}\oplus {\bf
R}i_k$, since there is the decomposition ${\sf A}= {\sf A}_0
i_0\oplus {\sf A}_1i_1\oplus ... {\sf A}_m i_m \oplus ...$ with
pairwise isomorphic real algebras ${\sf A}_0, {\sf A}_1,..., {\sf
A}_m,...$.
\par Let ${\cal B}(Y,{\cal A}_v)$ denote the algebra of all
bounded Borel functions from a topological space $Y$ into the
Cayley-Dickson algebra ${\cal A}_v$, let also ${\cal B}_u(Y,{\cal
A}_v)$ denote the algebra of all Borel functions from $Y$ into
${\cal A}_v$ with point-wise addition and multiplication of
functions and multiplication of functions $f$ on the left and on the
right on Cayley-Dickson numbers $a, b \in {\cal A}_v$, $~2\le v$.

\par {\bf 3. Lemma.} {\it If $T$ is a closed symmetrical operator in
a Hilbert space $X$ over the Cayley-Dickson algebra ${\cal A}_v$,
$~2\le v$, then ranges ${\cal R}(T\pm MI)$ of $(T\pm MI)$ are closed
for each $M\in {\cal S}_v$, when $v\le 3$ or $M\in \{ i_1,i_2,... \}
$. If $T$ is closed and $0\le ~ <Tz;z>$ for each $z\in {\cal D}(T)$,
then $T+I$ has a closed range.}
\par {\bf Proof.} Let $\mbox{}_nx$ be a sequence in ${\cal D}(T)$ so
that $ \{ (T\pm MI)\mbox{ }_nx : ~n \} $ tends to a vector $y\in X$.
But $<Tz;z>$ is real for each $z\in {\cal D}(T)$, hence $$\| z \| ^2
\le (<Tz;z>^2 + <z;z>^2)^{1/2} = |<(T\pm MI)z;z>|$$  $$\le \| (T\pm
MI)z \| \| z \| ,$$ consequently, $ \| \mbox{}_nx - \mbox{}_mx \|
\le \| (T\pm MI)(\mbox{}_nx-\mbox{ }_mx) \| $ and $\mbox{}_nx$
converges to some vector $x\in X$. On the other hand, the sequence
$\{ T\mbox{}_nx : ~ n \} $ converges to $\mp Mx+y$ and the operator
$T$ is closed, consequently, $x\in {\cal D}(T)$ and $Tx=\mp Mx+y$.
Indeed, $$M^*(Mx) = M^* \sum_jx_j(Mi_j)=\sum_jx_ji_j=x,$$ since $ \|
Mx \| = \| x \| $ and $M(x_ji_j) = x_j(Mi_j)$ and $M^*(Mi_j)=i_j$
for each $x_j\in X_j$ and $M\in {\cal S}_v$ for $2\le v \le 3$ or
$M\in \{ i_1, i_2,... \} $. Thus $(T\pm MI)x=y$, hence the operators
$(T\pm MI)$ have closed ranges for such $M$.
\par If an operator $T$ is closed and $0\le ~ <Tz;z>$ for each $z\in
{\cal D}(T)$, then \par $ \| z \| ^2 \le ~ <z;z> + <Tz;z> \le \|
(T+I)z\| \| z \| $ \\ and analogously to the proof above we get that
$(T+I)$ has a closed range.

\par {\bf 4. Proposition.} {\it If $T$ is a closed symmetrical operator
on a Hilbert space $X$ over the Cayley-Dickson algebra ${\cal A}_v$,
$~2 \le v$, then the following statements are equivalent:
\par $(1)$ an operator $T$ is self-adjoint;
\par $(2)$ operators $(T^*\pm MI)$ have $ \{ 0 \} $ as null space for each
$M\in {\cal S}_v$, when $v\le 3$ or $M\in \{ i_1,i_2,... \} $;
\par $(3)$ operators $(T\pm MI)$ have $X$ as range for every $M\in {\cal S}_v$,
when $v\le 3$ or $M\in \{ i_1,i_2,... \} $;
\par $(4)$ operators $(T\pm MI)$ have ranges dense in $X$ for all
$M\in {\cal S}_v$ with $v\le 3$ or $M\in \{ i_1,i_2,... \} $.}
\par {\bf Proof.} $(1)\Rightarrow (2)$. We have $<Tx;x> = <x;Tx> \in
\bf R$ for each $x\in {\cal D}(T)$, when $T^*=T$, hence
\par $<(T^*\pm MI)x;x> = <(T\pm MI)x;x> = <Tx;x> \pm M \| x \| ^2$,
\\
consequently, $<(T^*\pm MI)x;x> =0$ only when $x=0$, since $ \| Mx
\| = \| x \| $ for each $M\in {\cal S}_v$ with $v\le 3$ or $M\in \{
i_1, i_2,... \} $. Thus each operator $(T^*\pm MI)$ has $ \{ 0 \} $
as null space.
\par $(2)\Rightarrow (3)$. Ranges ${\cal R}(T\pm MI)$ are closed due
to Lemma 3. Therefore, it is sufficient to show that these ranges
are dense in a Hilbert space $X$ over the Cayley-Dickson algebra
${\cal A}_v$. But $<Tx;y> = \mp M <x;y>$, when $<(T\pm MI)x;y> =0$
for all $x\in {\cal D}(T)$, consequently, $y\in {\cal D}(T^*)$ and
$T^*y=\pm My$. Therefore, $y=0$, since the operators $T^*\pm MI$
have $ \{ 0 \} $ as null space. Thus the operators $(T\pm MI)$ have
dense ranges for each $M\in {\cal S}_v$ when $v\le 3$ or $M\in \{
i_1, i_2, ... \} $.
\par $(3)\Leftrightarrow (4)$. This follows from the preceding
demonstrations.
\par $(3)\Rightarrow (1)$. When $T$ is a closed and symmetrical
operator, one has that $T\subseteq T^*$ and a graph $\Gamma (T)$ is
a closed subspace of the closed $\bf R$-linear space $\Gamma (T^*)$.
The equality \par $<y;x> + <T^*y;Tx> =0$ \\ is valid for each $x\in
{\cal D}(T)$, when $(y,T^*y)\in \Gamma (T^*)$ is orthogonal to
$\Gamma (T)$. Operators $(T\pm MI)$ have range $X$, consequently,
there exists a vector $x\in {\cal D}(T)$ so that $Tx\in {\cal D}(T)$
and $y= (T+MI)(T-MI) x = (T^2+I)x$, since $T\subseteq T^*$ and
$T(MI)\subseteq (MI)T^*$. For such vector $x$ one gets
\par $<y;y> = <y;(T^2+I)x> = <y;x> + <T^*y;Tx>=0$, \\ hence
$<y;T^*y> = (0,0)$ and $\Gamma (T) = \Gamma (T^*)$ and hence
$T^*=T$. That is, the operator $T$ is self-adjoint.
\par {\bf 5. Note.} If an operator $T$ in a Hilbert space $X$ over
the Cayley-Dickson algebra ${\cal A}_v$, $~2\le v$, is self-adjoint,
the fact that operators $(T\pm MI)$ have dense everywhere defined
bounded inverses with bound not exceeding one follows from
Proposition 4 and the inequality at the beginning of the
demonstration of Lemma 3 for $M\in {\cal S}_v$ with $v\le 3$ or
$M\in \{ i_1, i_2,... \} $.
\par For an operator $T$ let $alg_{{\cal A}_v}(I,T)=:\sf Q$ be a family of
operators generated by $I$ and $T$ over the Cayley-Dickson algebra
${\cal A}_v$. Consider this family of operators on a common domain
${\cal D}^{\infty }(T) := \cap_{n=1}^{\infty } {\cal D}(T^n)$. Then
the family $\sf Q$ on ${\cal D}^{\infty }(T)$ can be considered as
an ${\cal A}_v$ vector space. Take the decomposition ${\sf Q}={\sf
Q}_0i_0\oplus {\sf Q}_1i_1\oplus ... \oplus {\sf Q}_mi_m\oplus ...$
of this ${\cal A}_v$ vector space with pairwise isomorphic real
vector spaces ${\sf Q}_0, {\sf Q}_1,...,{\sf Q}_m,... $ and for each
operator $B\in \sf Q$ put
$$(1)\quad B=\sum_j \mbox{}^jB\mbox{  with  }\mbox{}^jB={\hat \pi }^j(B)\in {\sf
Q}_ji_j$$ for each $j$, where ${\hat \pi }^j: {\sf Q}\to {\sf
Q}_ji_j$ is the natural $\bf R$ linear projection, real linear
spaces ${\sf Q}_mi_m$ and $i_m{\sf Q}_m$ are considered as
isomorphic. Evidently, there is the inclusion of domains of these
$\bf R$ linear operators ${\cal D}(\mbox{}^kB)\supset {\cal D}(B)$
for each $k$, particularly, ${\cal D}(\mbox{}^kT)\supset {\cal
D}(T)$.

\par {\bf 6. Lemma.} {\it Let $\{ \mbox{}_nE: ~ n \} $ be an
increasing sequence of ${\cal A}_v$ graded projections on a Hilbert
space $X$ over the Cayley-Dickson algebra ${\cal A}_v$ and let also
$G$ be an $\bf R$ homogeneous ${\cal A}_r$ additive operator with
dense domain $\bigcup_n \mbox{ }_nE(X)=: \cal E$ such that $G\mbox{
}_nE$ is a bounded self-adjoint operator on $X$, where $2\le v$.
Then $G$ is pre-closed and its closure $T$ is self-adjoint.
Moreover, if an operator $T$ is closed with core $\cal E$ and
$T\mbox{ }_nE$ is a bounded self-adjoint operator for each $n\in \bf
N$, then $T$ is self-adjoint.}
\par {\bf Proof.} For each $x, y \in \cal E$ there exists a natural
number $m$ so that \par $<Gx;y> = <G\mbox{ }_mEx;y> = <x;G\mbox{
}_mEy>=<x;Gy>$, \\ hence $y\in {\cal D}(G^*)$ and $G^*$ is densely
defined so that $G$ is pre-closed. \par Consider now the closure $T$
of $G$. For each $x, y \in {\cal D}(T)$ there exist sequences
$\mbox{}_nx$ and $\mbox{}_ny$ in $\cal E$ for which $\lim _n \mbox{
}_nx=x$, $ ~ \lim_n \mbox{ }_ny=y$, $~ \lim_n T \mbox{ }_nx=Tx$ and
$~ \lim_n T \mbox{ }_ny=Ty$. Therefore, the equalities
\par $<T\mbox{ }_nx; \mbox{ }_ny> = <T\mbox{ }_mE\mbox{ }_nx;
\mbox{ }_ny> = <\mbox{ }_nx; T\mbox{ }_mE\mbox{ }_ny> = <\mbox{
}_nx; T \mbox{ }_ny>$ \\ are satisfied, consequently, \par $\lim_n
<T\mbox{ }_nx; \mbox{ }_ny> = <Tx;y>$ and
\par $\lim_n <\mbox{ }_nx; T\mbox{ }_ny> = <x;Ty>$. \\
Thus one gets $<Tx;y> = <x;Ty>$ and that the operator $T$ is
symmetric. Mention that $(T\pm MI)\mbox{ }_nE(X) = \mbox{}_nE(X)$
for each $M\in {\cal S}_v$ with $v\le 3$ or $M\in \{ i_1, i_2,... \}
$, since $T\mbox{ }_nE$ is bounded and self-adjoint, for which the
operator $T\mbox{ }_nE \pm M\mbox{ }_nE$ has a bounded inverse on
$\mbox{}_nE(X)$. This implies that $T\pm MI$ has a dense range and
it coincides with $X$, consequently, $T$ is self-adjoint due to
Lemma 3 and Proposition 4.

\par {\bf 7. Theorem.} {\it Let $\mu $ be a $\sigma $-finite measure
$\mu : {\cal F} \to [0,\infty ]$ on a $\sigma $-algebra $\cal F$ of
subsets of a set $S$ and let $L^2(S,{\cal F},\mu ,{\cal A}_v)$ be a
Hilbert space completion of the set of all step $\mu $ measurable
functions $f: S\to {\cal A}_v$ with the ${\cal A}_v$ valued scalar
product
$$<f;g> =\int f(x){\tilde g}(x)\mu (dx)$$ for each $f, g \in
L^2(S,{\cal F},\mu ,{\cal A}_v)$, where $2\le v$. Suppose that $\sf
A$ is its left multiplication algebra $M_gf=gf$. Then an $\bf
R$-linear ${\cal A}_v$-additive operator $T$ is affiliated with $\sf
A$ if and only if a measurable finite $\mu $ almost everywhere on
$S$ function $g: S\to {\cal A}_v$ exists so that $T=M_g$ and
$g_k(t)f_j(t)=(-1)^{\kappa (j,k)} f_j(t)g_k(t)$ for $\mu $ almost
all $t\in S$ and each $f\in L^2(S,{\cal F},\mu ,{\cal A}_v)$ for
each $j, k=0,1,...$, where $\kappa (j,k)=0$ for $j=0$ or $k=0$ or
$j=k$, while $\kappa (j,k)=1$ for each $j\ne k\ge 1$. Moreover, an
operator $T\eta \sf A$ is self-adjoint if and only if $g$ is
real-valued $\mu $ almost everywhere on $S$.}
\par {\bf Proof.} If $g: S\to {\cal A}_v$ is a $\mu $ measurable
$\mu $ essentially bounded on $S$ function, then $M_g$ is a bounded
$\bf R$-linear ${\cal A}_v$-additive operator on $L^2(S,{\cal F},\mu
,{\cal A}_v)$, since \par $\| M_gf \|_2 \le \| g \| _{\infty } \| f
\|_2$, where \par $\| g \| _{\infty } = ess ~\sup_{x\in S} |g(x)|$
and $ \| f \|_2 := \sqrt{<f;f>}$. That is $M_g\in \sf A$.
\par Each operator $G\in \sf A$ is an $({\cal A}_v)_{{\bf C}_{\bf i}}$
combination of unitary operators in $\sf A$ (see Theorem II.2.20
\cite{ludopalglamb}) and $UT\subseteq TU$ for each unitary operator
$U\in \sf A$, since $T\eta \sf A$ and ${\sf A}\subseteq {\sf
A}^{\star }$, where ${\bf C}_{\bf i} = {\bf R}\oplus {\bf R}{\bf
i}$. If $F$ is an ${\cal A}_v$ graded projection operator
corresponding to the characteristic function $\chi _P$ of a $\mu
$-measurable subset $P$ in $S$, this implies the inclusion
$FT\subseteq TF$, consequently, $Ff\in {\cal D}(T)$ for each $f\in
{\cal D}(T)$. If $f\in {\cal D}(T)$ and $\mbox{}_nF$ corresponds to
the multiplication by the characteristic function $\chi _K$ of the
set $K = \{ x\in S: ~ |f(x)|\le n \} $, then $\mbox{}_nF$ is an
ascending sequence of projections in the algebra $\sf A$ converging
to the unit operator $I$ relative to the strong operator topology,
since $f$ is finite almost everywhere, $ \mu \{ x: ~ | f(x)| =\infty
\} =0$. Therefore, $\mbox{}_nFf\in \cal E$ for each $n$, where $\cal
E$ denotes the set of all $\mu $ essentially bounded functions $f\in
{\cal D}(T)$. Moreover, the limits exist:
$$\lim_n \mbox{ }_nFf=f\mbox{ and } \lim_n T\mbox{ }_nFf= \lim_n
\mbox{ }_nFTf = Tf.$$  Thus ${\cal E} = \bigcup_n \mbox{ }_nF{\cal
D}(T)$, where $\bigcup_n \mbox{ }_nF{\cal D}(T)$ is dense in ${\cal
D}(T)$.  Thus $\cal E$ is a core of $T$.
\par Each step function $u: S\to {\cal A}_v$ on $(S,{\cal F})$
has the form $$u(s) = \sum_{l=1}^m c_l \chi _{B_l},$$  where $B_l\in
\cal F$ and $c_l\in {\cal A}_v$ for each $l=1,...,m$, $~m\in \bf N$,
where $\chi _B$ denotes the characteristic function of a subset $B$
in $S$ so that $\chi _B(s)=1$ for each $s\in B$ and $\chi _B(s)=0$
for all $s$ outside $B$. A subset $N$ in $S$ is called $\mu $ null
if there exists $H\in \cal F$ so that $N\subset H$ and $\mu (H)=0$.
If consider an algebra ${\cal F}_{\mu }$ of subsets in $S$ which is
the completion of $\cal F$ by $\mu $ null subsets, then each step
function in $L^2(S,{\cal F},\mu ,{\cal A}_v)$ may have $B_l\in {\cal
F}_{\mu }$ for each $l$.
\par For each functions $f, g \in \cal E$ there are the equalities
\par $(1)$ $((f_ji_j)\mbox{ }^kT)g_k= (M_{f_ji_j}\mbox{ }^kT)g_k
=(-1)^{\kappa
(j,k)} (\mbox{}^kT M_{f_ji_j})g_k = \pm \mbox{}^lT f_jg_k$\par $ =
(-1)^{\kappa (j,k)} M_{g_ki_k} \mbox{ }^jTf_j =
(-1)^{\kappa (j,k)}g_ki_k [\mbox{ }^jTf_j]$, \\
where $i_ji_k=\pm i_l$, $ ~ \kappa (j,k)=0$ for $j=0$ or $k=0$ or
$j=k$, while $\kappa (j,k)=1$ for each $j\ne k\ge 1$, $$f=\sum_j
f_ji_j$$ with real-valued functions $f_j$ for each $j$. Let $S_k\in
\cal F$ be a sequence of pairwise disjoint subsets with $0<\mu
(S_k)<\infty $ for each $k$ and with union $\bigcup_k S_k=S$. For
the characteristic function $\chi _{S_k}$ a sequence $ \{ f^{k,j}: ~
j\in {\bf N} \} \subset \cal E$ of real-valued functions exists
converging to $\chi _{S_k}$ in $L^2(S,{\cal F},\mu ,{\cal A}_v)$.
The set $S_k^0 := \{ t\in S_k: ~ f^{k,j}(s)=0 ~ \forall j \} $ has
$\mu $ measure zero, since $$0=\lim_j M_gf^{k,j} =M_g\chi
_{S_k}=g,\mbox{  where  }g=\chi _{S_k^0}.$$ Put $h(s) =
[(Tf^{k,j})(s)][f^{k,j}(s)]^{-1}$ for $s\in S_k\setminus S_k^0$,
where $j$ is the least natural number so that $f^{k,j}(s)\ne 0$.
Thus $h$ is a measurable function defined $\mu $ almost everywhere
on $S$. \par In accordance with Formula $(1)$ the equality
$$(2)\quad f^{k,j}(s)[Tf]_l(s)=\sum_{p,q; ~ i_qi_p=i_l} \{
(-1)^{\kappa (q,p)} f_p(s)i_p[Tf^{k,j}]_q(s)i_q+
f_q(s)i_q[Tf^{k,j}]_p(s)i_p \} $$ is accomplished for all $k, j\in
\bf N$ and $l =0, 1, 2,...$ except for a set of measure zero,
consequently, $[Tf](s) = h(s)f(s)=M_hf(s)$ almost everywhere on $S$
so that $\mbox{}^kTf=[Tf]_ki_k$ for each $k$.
\par The operator $M_h$ is closed and affiliated with $\sf A$ as
follows from
the demonstration above. On the other hand, $M_h$ is an extension of
$T|_{\cal E}$, consequently, $T\subseteq M_h$. \par The family of
all functions $z\in L^2(S,{\cal F},\mu ,{\cal A}_v)$ vanishing on $
\{ s\in S: ~ |h(s)|>m \} $ for some $m=m(z)\in \bf N$ forms a core
for $M_h$. For such a function $z$ take a sequence $f^k$ of
functions in $\cal E$ tending to $z$ in $L^2(S,{\cal F},\mu ,{\cal
A}_v)$. Evidently $f^k$ can be replaced by $yf^k$, where $y$ is the
characteristic function of the set $ \{ s\in S: ~ z(s)\ne 0 \} $ of
all points at which $z$ does not vanish. Thus we can choose a
sequence $f^k$ vanishing for each $k$ when $z$ does. Therefore,
\par $(3)$ $\lim_k Tf^k = \lim_k M_hM_yf^k = M_hM_yz=M_hz$, \\ since
the operator $M_hM_y$ is bounded.
\par For a closed operator $T$ this implies that $z\in {\cal D}(T)$
and $Tz=M_hz$, consequently, $T=M_h$.
\par If an operator $T$ is self-adjoint, then $M_{yh}$ is a bounded
self-adjoint operator, hence the function $yh$ is real-valued $\mu $
almost everywhere on $S$.
\par For a bounded multiplication operator $M_g$ the ${\cal A}_v$
graded projections $\mbox{}_tE$ corresponding to multiplication by
characteristic function of the set $\{ s\in S: ~g(s)\le t \} $ forms
a spectral ${\cal A}_v$ graded spectral resolution $\{ \mbox{}_tE: ~
t \in {\bf R} \} $ of the identity for the operator $T$ (see also
Theorem 2.28 \cite{lujmsalop}).

\par {\bf 8. Definition.} Suppose that $V$ is an extremely
disconnected compact Hausdorff topological space and $V\setminus W$
is an open dense subset in $V$. If a function $f: V\setminus W \to
{\cal A}_v$ is continuous and $$\lim_{x\to y} |f(x)|=\infty $$ for
each $y\in W$, where $x\in V\setminus W$, $~ 1\le v$, then $f$ is a
called a normal function on $V$. \par If a normal function on $V$ is
real-valued it will be called a self-adjoint function on $V$.
\par The families of all normal and self-adjoint functions on $V$ we
denote by ${\cal N}(V,{\cal A}_v)$ and ${\cal Q}(V)$ respectively,
let also
$$W_+ :=W_+(f) := \{ y\in W: ~ \lim_{x\to y} f(x)=\infty \} ,$$ $$W_-
:=W_-(f) := \{ y\in W: ~ \lim_{x\to y} f(x)= - \infty \} $$ for a
self-adjoint function $f$ on $V$.

\par {\bf 9. Lemma.} {\it Let $f$ and $g$ be normal functions on $V$
(see Definition 8) defined on $V\setminus W_f$ and $V\setminus W_g$
respectively so that $f(x)=g(x)$ for each $x$ in a dense subset $U$
in $V\setminus (W_f\cup W_g)$. Then $W_f=W_g$ and $f=g$.}
\par {\bf Proof.} The subset $V\setminus (W_f\cup W_g)$ is dense
in $V$, consequently, $U$ is dense in $V$. If $y\in W_g$, then for
each $N>0$ there exists a neighborhood $E$ of $y$ in $V$ so that
$|g(x)|>N$ for each $x\in E\cap (V\setminus W_g)$, hence
$|f(x)|=|g(x)|>N$ for each $x\in E\cap (V\setminus (W_g\cup W_f))$,
consequently, $W_g\subset W_f$ and symmetrically $W_f\subset W_g$.
Thus $W_f=W_g$ and $f-g$ is defined and continuous on $V\setminus
W_f$ and is zero on $U$, hence $f=g$.
\par {\bf 10. Lemma.} {\it Suppose that $T$ is a self-adjoint operator
acting on a Hilbert space $X$ over either the quaternion skew field
or the octonion algebra ${\cal A}_v$, $~2\le v\le 3$, so that $T$ is
affiliated with some quasi-commutative von Neumann algebra ${\sf A}$
over ${\cal A}_v$, where $\sf A$ is isomorphic to $C(\Lambda ,{\cal
A}_v)$ with an extremely disconnected compact Hausdorff topological
space $\Lambda $. Then there is a unique self-adjoint function $h$
on $\Lambda $ such that $h\hat \cdot e\in C(\Lambda ,{\cal A}_v)$
and a function $h\hat \cdot e$ represents $TE$, when $E$ is an
${\cal A}_v$ graded projection for which $TE\in L_q(X)$ is a bounded
operator on $X$, a function $e\in C(\Lambda ,{\cal A}_v)$
corresponds to $E$, $~ h\hat \cdot e(x)=h(x)$ for $e(x)=1$, while $~
h\hat \cdot e(x)=0$ otherwise. There exists an ${\cal A}_v$ graded
resolution of the identity $ \{ \mbox{}_bE: ~ b \} $ so that
$\bigcup _{n=1}^{\infty } \mbox{}_nF(X)$ is a core for $T$, where
$\mbox{}_nF := \mbox{}_nE-\mbox{ }_{-n}E$ and \par $(1)$
$Tx=\int_{-n}^n d\mbox{}_bE.b x$ \\ for every $x\in \mbox{ }_nF(X)$
and each $n$ in the sense of norm convergence of Riemann sums.}
\par {\bf Proof.} Take $Y=X\oplus {\bf i}X$, where $\bf i$ is a generator
commuting with $i_j$ for each $j$ such that ${\bf i}^2=-1$. Take an
extension $T$ onto ${\cal D}(T)\oplus {\cal D}(T){\bf i}$ so that
$T(x+y{\bf i})=Tx+(Ty){\bf i}$ for each $x, y \in {\cal D}(T)$. Then
${\cal R}(T+{\bf i}I)=Y$ and ${\cal R}(T-{\bf i}I)=Y$, where ${\cal
R}(T\pm {\bf i}I)=(T\pm {\bf i}I)Y$, $~ker (T\pm {\bf i}I)= \{ 0 \}
$ and inverses $B_{\pm } := (T\pm {\bf i}I)^{-1}$ are everywhere
defined on $Y$ and of norm not exceeding one in accordance with \S
II.2.74 and Proposition II.2.75 \cite{ludopalglamb}. Then the
equalities \par $<B_+(T+{\bf i}I)x; (T-{\bf i}I)y> = <x;(T-{\bf
i}I)y> = <(T+{\bf i}I)x;y> = <(T+{\bf i}I)x,B_-(T-{\bf i}I)y>$ \\
are satisfied for each $x$ and $y\in {\cal D}(T)$, since the
operator $T$ is self-adjoint, hence $B_- = B_+^*$.
\par Then an arbitrary vector $z\in X$ has the form $z=B_- B_+ x$ for
the corresponding vector $x\in {\cal D}(T)$ with $Tx\in {\cal
D}(T)$, since ${\cal R}(B_{\pm }) = Y$. In the latter case we get
$B_- B_+ x = (TT- {\bf i}IT+T{\bf i}I + I)x = B_+ B_- x$, since to
$\mbox{}_nF$ a real-valued function $\mbox{}_nf\in C(\Lambda ,{\bf
R})$ corresponds for each $n$. On the other hand, $B_-=B_+^*$,
consequently, the operator $B_+$ is normal.
\par Consider a quasi-commutative von Neumann algebra $\breve{\sf A}$
over the complexified algebra $({\cal A}_v)_{\bf C_i}$ containing
$I$, $B_+$ and $B_- \in \sf A$. If $U$ is a unitary operator in
${\sf A}^{\star }$ (see \S II.2.71 \cite{ludopalglamb}), then
$$Ux=UB_+(T+{\bf i}I)x=
(B_+) U (T+{\bf i}I)x\mbox{ hence}$$
$$(T+{\bf i}I)Ux= U(T+{\bf i}I)x,$$ consequently, $T$ is affiliated with
$\breve{\sf A}$ by Definition 1, since $({\bf i}I)Ux= U({\bf i}I)x$
for a unitary operator $U$ as follows from the definition of a
unitary operator and the ${\cal A}_v$ valued scalar product on $X$
extended to the $({\cal A}_v)_{\bf C_i}$ valued scalar product on
$Y$:
\par $<a+b{\bf i};c+q{\bf i}> = (<a;c> + <b;q>) +(<b;c> - <a;q>){\bf
i}$
\\ for each vectors $a, b, c, q\in X$.
\par The algebra $\breve{\sf A}$ has the decomposition $\breve{\sf A}
= {\sf A}^0\oplus {\sf A}^1{\bf i}$, where ${\sf A}^0$ and ${\sf
A}^1$ are quasi-commutative isomorphic algebras over either the
quaternion skew field or the octonion algebra ${\cal A}_v$ with
$2\le v \le 3$.
\par In view of Theorem 2.24 \cite{lujmsalop} $\breve{\sf A}$ is isomorphic with
$C(\Lambda ,({\cal A}_v)_{\bf C_i})$ for some extremely disconnected
compact Hausdorff topological space $\Lambda $, since the generator
$\bf i$ has the real matrix representation ${\bf i} ={{0 ~~ 1
}\choose {-1 ~ 0}}$ and the real field is the center of the algebra
${\cal A}_v$ with $2\le v\le 3$. Let $f_+$ and $f_-$ be functions
corresponding to $B_+$ and $B_-$ respectively. Put $h_{\pm } =
1/f_{\pm }$ at those points where $f_{\pm }$ does not vanish.
Therefore, these functions $h_{\pm }$ are continuous on their
domains of definitions as well as $h=\frac{1}{2}(h_+ + h_-)$. The
function $h$ corresponds to $T$.
\par It will be demonstrated below that $h$ is real-valued and then
the ${\cal A}_v$ graded spectral resolution of $T$ will be
constructed.
\par At first it is easy to mention that $f_+^*=f_-$, since
$B_+^*=B_-$, consequently, $h$ is real-valued. The functions $f_+$
and $f_-$ are continuous and conjugated, hence the set $W :=
f_+^{-1}(0)=f_-^{-1}(0)$ is closed. If the set $W $ contains some
non-void open subset $J$, then $cl (J)\subset W$ so that $cl (J)$ is
clopen (i.e. closed and open simultaneously) in $\Lambda $,
consequently, $W$ is nowhere dense in $\Lambda $. The projection $P$
corresponding to this subset $cl (J)$ would have the zero product
$B_+P=0$ contradicting the fact that $ker (B_+)= \{ 0 \} $. Thus
each point $x\in W$ is a limit point of points $y\in \Lambda
\setminus W$, but $h_+$ and $h_-$ are defined on the latter set.
Therefore, the set $\Lambda \setminus W$ is dense in $\Lambda $ and
the functions $h_+$ and $h_-$ are defined on $\Lambda \setminus W$.
\par Then $TB_+B_- y = (T+{\bf i}I-{\bf i}I) B_+B_-y = B_-y - ({\bf
i}I)B_+B_-y$ for each $y\in \Lambda $, hence $TB_+B_- = B_- - ({\bf
i}I)B_+B_-$ and analogously $TB_-B_+ = B_+ + ({\bf i}I)B_-B_+$ and
inevitably \par $(2)$ $2({\bf i}I)B_+B_- = B_- - B_+$ and \par $(3)$
$TB_+B_- = \frac{1}{2} (B_+ + B_-)$.
\par From Formula $(2)$ and the definitions of functions $f_{\pm }$ and $h$
one gets $(h(y)+{\bf i})^{-1} = f_+(y)$ and  $(h(y)-{\bf i})^{-1} =
f_-(y)$ for each $y\in \Lambda \setminus W$, since $(h(y)+{\bf
i})f_+(y) = \frac{1}{2} +\frac{1}{2f_-}f_+ +{\bf i}f_+ =
\frac{1}{2}+\frac{w-{\bf i}}{2}\frac{1}{w+{\bf i}} +{\bf
i}\frac{1}{w+{\bf i}}= 1 -\frac{2{\bf i}}{2}\frac{1}{w+{\bf i}}+{\bf
i}\frac{1}{w+{\bf i}}=1$, where $w\pm {\bf i}=\frac{1}{f_{\pm }}$
corresponds to $(T\pm {\bf i}I)$.  But $f_+(y)$ tends to zero when
$y\in \Lambda \setminus W$ tends to $x\in W$, consequently,
$\lim_{y\to x} |h(y)|=\infty $. This means that $h$ is a
self-adjoint function on $\Lambda $.
\par We put $U_b := \{ x\in \Lambda \setminus W: ~h(x)>b \} $ and
$V_b=U_b\cup W_+(h)$ for $b\in {\bf R}$ (see Definition 8). The
function $h$ is continuous on the open subset $\Lambda \setminus W$
in $\Lambda $, hence the subset $U_b$ is open in $\Lambda $. For a
point $x\in W_+$ there exists an open set $Q$ containing $x$ such
that $h(y)>\max (b,0)$ for each $y\in Q\cap (\Lambda \setminus W)$,
since $\lim_{y\to x} |h(y)|=\infty $, where $W_+=W_+(h)$ (see
Definition 8). Therefore, this implies the inclusion $Q\cap (\Lambda
\setminus W)\subset U_b\subset V_b$.
\par Suppose that there would exist a point $z\in Q\cap W_-$, where $W_- =
W_-(h)$, then a point $y\in Q$ with $h(y)<0$ would exist
contradicting the choice of $Q$. This means that $Q\cap W\subset
W_+$ and $Q\subset V_b$, consequently, $V_b$ is open in $\Lambda $.
\par We consider next the subset $\Lambda _b := \Lambda \setminus cl
(V_b)$. In accordance with \S 2.24 \cite{lujmsalop} the set $\Lambda
_b$ contains every clopen subset $K_b := \{ y\in \Lambda : ~ h(y)\le
b \} $, where $h(y)= - \infty $ for $y\in W_-$ so that $W_-\subset
K_b$ for each $~b\in \bf R$. To demonstrate this suppose that $y\in
V_b$, then $y\in \Lambda \setminus K_b$ and $cl (V_b)\subset \Lambda
\setminus K_b$ and $K_b\subset \Lambda \setminus cl (V_b)=\Lambda
_b$, since $\Lambda \setminus K_b$ is closed. Moreover, we have
$h(y)\le b$ for each $y\in \Lambda _b\cap (\Lambda \setminus W)$,
since $y\notin U_b$, while the set $~\Lambda _b$ is clopen in
$\Lambda $. Therefore, $y\in W\setminus W_+$ and $h(y)\le b$ for
each $y\in \Lambda _b\cap W$ and $W_-=W\setminus W_+$. Thus $\Lambda
_b$ is the largest clopen subset in $\Lambda $ so that $h(y)\le b$
and $\Lambda _b\cap W=W_-$.
\par Denote by $e_b$ the characteristic function of the subset $\Lambda _b$
and $\mbox{}_bE$ be an ${\cal A}_v$ graded projection operator in
$\sf A$ corresponding to $e_b$, where $b\in \bf R$ (see also \S 2.24
\cite{lujmsalop}). The subset $W$ is nowhere dense in $\Lambda $,
consequently, $\vee_b e_b=1$ and $\wedge _b e_b =0$ such that
$\vee_b \mbox{ }_bE =1$ and $\wedge _b \mbox{ }_bE =0$. Thus $\{
\mbox{}_bE: ~ b \} $ is the ${\cal A}_v$ graded resolution of the
identity so that $\mbox{}_bE_s=(sI)\mbox{ }_bE=\mbox{}_bE(sI)$
corresponds to $e_bs=se_b$ for each marked Cayley-Dickson number
$s\in {\cal A}_v$. This resolution of the identity is unbounded when
$h\notin C(\Lambda ,{\cal A}_v)$.
\par Put $F = \mbox{}_bE - \mbox{ }_aE$ for $a<b\in \bf R$,
hence $e_b - e_a=:u$ is the characteristic function of $\Lambda
_b\setminus \Lambda _a$ corresponding to $F$. Therefore, the
inclusion $\Lambda _b\setminus \Lambda _a\subset \Lambda \setminus
W$ follows and $f_+(y)f_-(y)\ne 0$ when $u(y)=1$, since $\Lambda
_b\cap W=\Lambda _a\cap W=W_-$. Then $$(4)\quad h(y)=
\frac{f_++f_-}{2f_+f_-}(y)$$ for each $y\in \Lambda \setminus W$.
The function $f_+f_-$ is continuous and vanishes nowhere on the
clopen subset $\Lambda _b\setminus \Lambda _a$, consequently, a
positive continuous function $\psi $ on $\Lambda $ exists so that
$\psi f_+f_-=u$ and $\psi u=\psi $. Consider an element $\Psi $ of
the algebra $\sf A$ corresponding to $\psi $, hence
\par $(5)$ $\Psi B_+B_-=F$.
\par On the other hand, from the construction above it follows that
$a\le h(y)\le b$ for each $y\in \Lambda _b\setminus \Lambda _a$ and
from Formula $(4)$ one gets
\par $af_+f_-u\le \frac{(f_+ + f_-)u}{2} \le b f_+f_-u$ and
\par $a\psi f_+f_-u=au\le \frac{(f_+ +f_-)\psi u}{2}= \frac{(f_+ +
f_-)\psi }{2} \le b\psi f_+f_-u=bu$, since the real field is the
center of the Cayley-Dickson algebra ${\cal A}_v$. Thus \par $(6)$
$aF\le \frac{(B_++B_-)\Psi }{2}\le bF$. \par Therefore, Formulas
$(3,5,6)$ imply that
\par $(7)$ $aF\le TF\le bF$.
\par Therefore, the operator $TF$ is bounded and the element
$h\hat \cdot u\in C(\Lambda ,{\cal A}_v)$ corresponds to it due to
$(3-5)$. If $E$ is an ${\cal A}_v$ graded projection belonging to
$\sf A$ so that $TE\in L_q(X)$ and $U$ is a unitary operator in
${\sf A}^{\star }$ such that $U^{-1}TEU=U^{-1}TUE=TE$ one has $TE\in
\sf A$ (see \S II.2.71 \cite{ludopalglamb}). Let each function
$\mbox{}_ng\in C(\Lambda ,{\cal A}_v)$ be corresponding to the
${\cal A}_v$ graded projection $\mbox{}_nF$ and $\Lambda _n :=
\mbox{}_ng^{-1}(1)$. Then $\bigcup_n \Lambda _n$ is dense in
$\Lambda $, since $\vee _n \mbox{ }_nF=I$. If a function $e\in
C(\Lambda ,{\cal A}_v)$ corresponds to $E$, then $(h\hat \cdot
\mbox{ }_ng)e$ corresponds to the operator $T\mbox{ }_nFE$, where
$\mbox{}_nFE= E\mbox{ }_nF$.
\par Suppose that $e(y)=1$ for some $y\in \Lambda $. For each (open)
neighborhood $H$ of $y$ there exist $n\in \bf N$ and $x\in \Lambda
_n$ such that $x\in H$, hence $((h\hat \cdot \mbox{ }_ng)e)(x)=h(x)$
and $|h(x)| \le \| TE\mbox{ }_nF \| \le \| TE \| $. Thus $y\notin W$
and $|h(y)|\le \| TE \| $ and hence $h\hat \cdot e \in C(\Lambda
,{\cal A}_v)$. Then one gets also $(h_j\hat \cdot \mbox{ }_ng_k)e_l
= h_j\hat \cdot (\mbox{ }_ng_ke_l)=(h_j\hat \cdot e_l)\mbox{ }_ng_k$
for each $j, k, l$. If $s\in C(\Lambda ,{\cal A}_v)$ represents
$TE$, then $s\mbox{ }_ng$ represents $TE\mbox{ }_nF$, that is
$s\mbox{ }_ng= (h\hat \cdot e)\mbox{ }_ng$ for every $n$ and
$s=h\hat \cdot e$.
\par On the other hand, $(2\mbox{ }_nF -I)T (2\mbox{ }_nF -I)=T$,
since $T$ and $\mbox{}_nF$ are $\bf R$-linear operators, also
$(\mbox{}_nF-I)(X)= (I-\mbox{}_nF)(X)$ and $X=\mbox{}_nF(X)\oplus
(\mbox{}_nF-I)(X)$, hence $\mbox{}_nFT\subset T\mbox{ }_nF$.
Moreover, $\lim_n \mbox{ }_nFx=x$ and \par $\lim_n T\mbox{ }_nFx =
\lim_n \mbox{ }_nFTx=T x$ \\
hence $\bigcup_n \mbox{ }_nF(X)$ is a core for $T$. \par In view of
Theorem I.3.9 \cite{ludopalglamb} we get that $ \{ E\mbox{
}_bE|_{E(X)} : ~ b\in {\bf R} \} $ is an ${\cal A}_v$ graded
resolution of the identity on $E(X)$, since $(h\hat \cdot e)ee_b\le
b e e_b$ and $b(e-ee_b) \le (h\hat \cdot e)(e-ee_b)$. Applying
Theorem I.3.6 \cite{ludopalglamb} for $\mbox{}_nF$ in place of $E$
for $T\mbox{ }_nF|_{\mbox{}_nF(X)}$ and $ \{ \mbox{}_nF \mbox{
}_bE|_{\mbox{}_nF(X)}: ~ b \} $ leads to the formula $$Tx =
\int_{-n}^n d\mbox{}_bE.b x$$ for each $x\in \mbox{}_nF(X)$ and
every $n\in \bf N$.
\par {\bf 11. Remark.} Lemma 10 means that
$$\lim_n \int_{-n}^n d\mbox{}_bE.b x = \lim_n\int_{-n}^n d\mbox{}_bE.b \mbox{ }_nFx =
\lim_nT\mbox{ }_nFx=Tx=\int_{-\infty }^{\infty }d\mbox{}_bE.b x$$
for each $x\in {\cal D}(T)$ interpreting the latter integral as
improper.
\par Under conditions imposed in Lemma 10 one says that the function $h\in {\cal Q}(\Lambda )$
represents an affiliated operator $T\eta \sf A$.
\par Mention that an ${\cal A}_v$ graded projection operator
$\mbox{}_bE$ is $\bf R$-homogeneous and ${\cal A}_v$-additive, hence
${\bf R}$-linear. Therefore, in the particular case of a real-valued
Borel function $h(b)$ one has $d\mbox{}_bE.h(b)x=h(b)d\mbox{}_bE.1x$
or simplifying notation $h(b)d\mbox{}_bEx$.
\par {\bf 12. Lemma.} {\it Suppose that $\sf A$ is a quasi-commutative von Neumann algebra
over ${\cal A}_v$, where $~2\le v\le 3$, so that $\sf A$ is
isomorphic to $C(\Lambda ,{\cal A}_v)$ for some extremely
disconnected compact Hausdorff topological space $\Lambda $. Then
each function $h\in {\cal Q}(\Lambda )$ represents some self-adjoint
operator $T$ affiliated with $\sf A$.}
\par {\bf Proof.} From \S 10 it follows that a self-adjoint function $h$
determines an ${\cal A}_v$ graded resolution of the identity $ \{
\mbox{}_bE: ~ b\in {\bf R} \} $ in $\sf A$. Moreover, $h\hat \cdot
g_n\in C(\Lambda ,{\cal A}_v)$, where $g_n=e_n - e_{-n}$, $~e_n\in
C(\Lambda ,{\cal A}_v)$ represents $\mbox{}_nE$. Consider an
operator $\mbox{}_nT$ corresponding to $h\hat \cdot g_n$. Certainly
one has $(h\hat \cdot g_m)g_n=h\hat \cdot g_n$ for each $n\le m$,
consequently, $\mbox{}_mT\mbox{}_nF=\mbox{}_nT$, where $\mbox{}_nF$
corresponds to $g_n$. Put $Gx=\mbox{}_nTx$ for every vector $x\in
\mbox{ }_nF(X)$ and $n\in \bf N$. Therefore, $G$ is an $\bf R$
linear ${\cal A}_v$ additive operator on $\bigcup_{n=1}^{\infty }
\mbox{}_nF(X)=:\cal K$. Therefore, the operator $G$ is pre-closed
and its closure $T$ is a self-adjoint operator with core $\cal K$ as
Lemma 6 asserts. For a unitary operator $U$ in ${\sf A}^{\star }$
and $x\in \mbox{ }_nF(X)$ we get $Ux\in \mbox{ }_nF(X)$, hence \par
$TUx = \mbox{ }_nTUx= U\mbox{ }_nTx $.
\\ Therefore, $T\eta \sf A$ due to Definition 1 and
Remark 2.
\par If $u\in {\cal Q}(\Lambda )$ represents $T$, then $u\hat \cdot
g_ n$ represents $T\mbox{ }_nF$ in accordance with Lemma 10,
consequently, $h\hat \cdot g_n = u\hat \cdot g_n$ for each $n$. In
view of Lemma 9 $h=u$, since $h$ and $u$ are consistent on a dense
subsets in $\Lambda $. Thus the function $h$ represents the
self-adjoint operator $T$.
\par {\bf 13. Lemma.} {\it Suppose that $\{ \mbox{}_bE: ~ b \} $
is an ${\cal A}_v$ graded resolution of the identity on a Hilbert
space $X$ over the Cayley-Dickson algebra ${\cal A}_v$, $~2\le v\le
3$, also $\sf A$ is a quasi-commutative von Neumann algebra over
${\cal A}_v$ so that $\sf A$ contains $\{ \mbox{}_bE: ~ b \} $,
where $2\le v \le 3$. Then there exists a self-adjoint operator $T$
affiliated with $\sf A$ so that
$$(1)\quad Tx = \int_{-n}^n d\mbox{}_bE.b x$$ for each $x\in \mbox{}_nF(X)$
and every $n\in \bf N$, where $\mbox{}_nF= \mbox{}_nE-\mbox{
}_{-n}E$, and $\{ \mbox{}_bE: ~ b \} $ is the ${\cal A}_v$ graded
resolution of the identity for $T$ given by Lemma 10.}
\par {\bf Proof.} Take a function $e_b\in C(\Lambda ,{\cal A}_v)$
corresponding to $\mbox{}_bE$ and a subset $\Lambda _b=e_b^{-1}(1)$
clopen in $\Lambda $. Consider the subsets $W_+ := \Lambda \setminus
\bigcup_b \Lambda _b$ and $W_- := \bigcap_b \Lambda _b$. From their
definition it follows that $W_+$ and $W_-$ are closed in $\Lambda $.
These two subsets $W_{\pm }$ are nowhere dense in $\Lambda $, since
$\wedge _b \mbox{ }_bE=0$ and $\vee _b \mbox{ }_bE=I$, hence $W :=
W_+\cup W_-$ is nowhere dense in $\Lambda $. For a point $x\in
\Lambda \setminus W$ we put $h(x) := \inf \{ b: ~ x\in \Lambda _b \}
$. Given a positive number $\epsilon >0$ and $y\in \Lambda \setminus
W$ so that $h(y)=b$ one gets $|h(z)-h(y)|\le \epsilon $ for each
$z\in \Lambda _{b+\epsilon }\setminus \Lambda _{b-\epsilon }$. This
means that the function $h$ is continuous on $\Lambda \setminus W$.
Then $\lim_{y\to x} h(y)=\pm \infty $ for $x\in W_{\pm }$, where
$y\in \Lambda \setminus W$. Thus the function $h$ is self-adjoint
$h\in {\cal Q}(\Lambda )$ and by Lemma 12 corresponds to a
self-adjoint operator $T$ affiliated with $\sf A$. Certainly we get
that $\{ \mbox{}_bE: ~ b \} $ is the ${\cal A}_v$ graded resolution
of the identity for the operator $T$ and Formula $(1)$ is valid,
since $\Lambda _b$ is the largest clopen subset in $\Lambda $ on
which the function $h$ takes values not exceeding $b$.
\par If $\Psi $ is another such clopen subset and $e$ is its
characteristic function, $E\in \sf A$ corresponds to $e$, then $\Psi
\subset \Lambda _c$ for each $c\ge b$. This leads to the conclusion
that $E\le \wedge _{c>b} \mbox{ }_cE$ and $\Psi \subset \Lambda _b$.
\par {\bf 14. Lemma.} {\it Suppose that $T$ is a closed operator on
a Hilbert space $X$ over either the quaternion skew field or the
octonion algebra ${\cal A}_v$ with $2\le v \le 3$ and $\{
\mbox{}_bE: ~ b \} $ is an ${\cal A}_v$ graded resolution of the
identity on $X$, where ${\cal E} := \bigcup_n \mbox{ }_nF(X)$ is a
core for $T$ while $\mbox{}_nF=\mbox{}_nE - \mbox{ }_{-n}E$ and
$$(1)\quad Tx = \int_{-n}^n d\mbox{}_bE.b  x$$ for each $x\in \mbox{}_nE(X)$
and all $n$, then $T$ is self-adjoint and $\mbox{}_bE$ is the ${\cal
A}_v$ graded resolution of the identity for $T$.}
\par {\bf Proof.} Formula $(1)$ implies that $T\mbox{ }_nF$ is
bounded and everywhere defined and is the strong operator limit of
finite real-linear combinations of $ \{ \mbox{}_bE: ~ b \} $. Thus
$\mbox{}_bE (T\mbox{ }_nF) =(T\mbox{ }_nF) \mbox{ }_bE$ and $T\mbox{
}_nF$ is self-adjoint, hence $T$ is self-adjoint by Lemma 6. For
each vector $x\in {\cal D}(T)$ there exists a subsequence $ \{ n_p:
~ p \} $ of natural numbers and a sequence $ \{ \mbox{}_p x : ~ p \}
$ of vectors such that \par $\lim_p \mbox{}_px = \lim_p
\mbox{}_{n_p}F \mbox{ }_px=x$ and $\lim_p T\mbox{ }_px = Tx$, \\
since $\cal E$ is a core for the operator $T$. Therefore,
$\mbox{}_nFT\subseteq T\mbox{ }_nF$ for each $n$, since \par
$\mbox{}_nFTx = \lim_p \mbox{}_nF (T\mbox{ }_{n_p}F) \mbox{ }_px =
\lim_p (T\mbox{ }_{n_p}F)\mbox{}_nF\mbox{ }_px = \lim_p
T\mbox{}_nF\mbox{ }_px = T\mbox{}_nFx$.
\par On the other hand, the limits exist:
\par $\lim_n ~ (T\mbox{ }_nF)\mbox{
}_bEx = \lim_n \mbox{ }_bE(T\mbox{ }_nF)x = \mbox{ }_bETx$ and
$\lim_n \mbox{ }_nF\mbox{ }_bEx =\mbox{ }_bEx$. The operator $T$ is
closed, hence $\mbox{}_bEx \in {\cal D}(T)$ and $T \mbox{ }_bEx =
\mbox{ }_bE Tx$. This leads to the conclusion that
$\mbox{}_bET\subseteq T\mbox{ }_bE$ and $(2\mbox{ }_bE-I)T(2\mbox{
}_bE-I)=T$, since these operators are $\bf R$ linear and
$X=\mbox{}_bE(X)\oplus (\mbox{}_bE-I)(X)$. Therefore,
\par $\mbox{}_bE(B_{\pm }) = (B_{\pm
}) \mbox{}_bE$ for each $b\in \bf R$.
\par Take the quasi-commutative von Neumann algebra
\par ${\sf G}={\sf G}(T) = cl [~ alg_{{\cal A}_v} \{ B_-, B_+; ~
\mbox{}_bE: ~ b \in {\bf R} \} \oplus alg_{{\cal A}_v} \{ B_-, B_+;
~ \mbox{}_bE: ~ b \in {\bf R} \}{\bf i} ]$. \\ In view of Lemma 10
the operator $T$ is affiliated with $\sf G$. But Lemma 13 means that
there exists a self-adjoint operator $H$ affiliated with $\sf G$ so
that $$Hx = \int_{-n}^n d\mbox{}_bE.b  x $$ for each $x\in \mbox{
}_nF(X)$ and every natural number $n$. Therefore, $H=T$ and $ \{
\mbox{}_bE: ~ b \} $ is the ${\cal A}_v$ graded resolution of the
identity for $T$, since $H|_{\cal E} = T|_{\cal E}$ and $\cal E$ is
a core for $T$ and $H$ simultaneously.
\par {\bf 15. Note.} The quasi-commutative von Neumann algebra
${\sf G}={\sf G}(T)$ generated by $B_-$ and $B_+$ in \S 14 with
which the self-adjoint operator $T$ is affiliated will be called the
von Neumann algebra generated by $T$.
\par {\bf 16. Theorem.} {\it Let $\sf A$ be a quasi-commutative
von Neumann algebra acting on a Hilbert space $X$ over either the
quaternion skew field or the octonion algebra ${\cal A}_v$, $~ 2\le
v\le 3$, let also ${\cal Q}({\sf A})$ be a family of all
self-adjoint operators affiliated with $\sf A$. Suppose that $\sf A$
is isomorphic to $C(\Lambda ,{\cal A}_v)$ for an extremely
disconnected compact Hausdorff topological space $\Lambda $ and
${\cal Q}(\Lambda )$ is the family of all self-adjoint functions on
$\Lambda $. Then
\par $(a)$ there exists a bijective mapping $\phi $ from ${\cal
Q}({\sf A})$ onto ${\cal Q}(\Lambda )$ which is an extension of the
isomorphism of $\sf A$ with $C(\Lambda ,{\cal A}_v)$ for which $\phi
(T)\hat \cdot e$ corresponds to $TE$ for each projection $E$ in $\sf
A$ with $TE\in \sf A$, where $e\in C(\Lambda ,{\cal A}_v)$
corresponds to $E$ and $((\phi (T)\hat \cdot e)(y)= (\phi (T))(y)$
for $e(y)=1$ while $(\phi (T)\hat \cdot e)(y)=0$ for $e(y)=0$; \par
$(b)$ an ${\cal A}_v$ graded resolution $\{ \mbox{}_bE: ~ b \} $ of
the identity exists in the quasi-commutative von Neumann subalgebra
$\sf G$ generated by an operator $T$ in ${\cal Q}({\sf A})$ so that
\par $$(1)\quad Tx = \int_{-n}^n d\mbox{}_bE.bx$$ for each $x\in
\mbox{}_nF(X)$ and every natural number $n$, where $\mbox{}_nF =
\mbox{}_nE- \mbox{ }_{-n}E$ and $\bigcup_n \mbox{ }_nF(X)=:\cal E$
is a core for $T$; \par $(c)$ if $\{ \mbox{}_bQ: ~ b\in {\bf R} \} $
is an ${\cal A}_v$ graded resolution of the identity on $X$ so that
\par $$(2)\quad Tx = \int_{-n}^n d\mbox{}_bQ.bx$$ for each
$x\in \mbox{}_nP(X)$ and every natural number $n$, where $\mbox{}_nP
= \mbox{}_nQ- \mbox{ }_{-n}Q$ and $\bigcup_n \mbox{ }_nP(X)$ is a
core for $T$, then $\mbox{}_bE=\mbox{}_bQ$ for all $b$; \par $(d)$
if $\{ \mbox{}_bE : ~b \in {\bf R} \} $ is an ${\cal A}_v$ graded
resolution of the identity in $\sf A$, then there exists an operator
$T\in {\cal Q}({\sf A})$ for which Formula $(1)$ is valid;
\par $(e)$ if a function $e_b \in C(\Lambda ,{\cal A}_v)$
corresponds to $\mbox{}_bE$ and $\Lambda _b=e_b^{-1}(1)$, then
$\Lambda _b$ is the largest clopen subset of $\Lambda $ on which
$\phi (T)$ takes values not exceeding $b$ in the extended sense.}
\par The latter theorem is the reformulation of the results obtained
above in this section.
\par {\bf 17. Definition.} A closed densely defined operator $T$ in
a Hilbert space $X$ over the Cayley-Dickson algebra ${\cal A}_v$
will be called normal when two self-adjoint operators $T^*T$ and
$TT^*$ are equal, where $2\le v$.
\par {\bf 18. Remark.} For an unbounded $\bf R$ homogeneous
${\cal A}_v$ additive operators the following properties are
satisfied:
\par $(1)$ if $A\subseteq B$ and $C\subseteq D$, then
$A+C \subseteq B+D$;
\par $(2)$ if $A\subseteq B$, then $CA\subseteq CB$ and $AC\subseteq
BC$;
\par $(3)$ $(A+B)C=AC+BC$ and $CA+CB\subseteq C(A+B)$.
\par The latter inclusion does not generally reduce to the equality.
For example, an operator $C$ may be densely defined, but not
everywhere, one can take $A=I$ and $B=-I$. This gives $C(A+B)=0$ on
$X$, but $CA+CB$ is zero only on ${\cal D}(C)$.
\par These rules imply that if $CA\subseteq AC$ for each $C$ in some
family $\cal Y$, then $TA\subseteq AT$ for each sum of products of
operators from $\cal Y$. Apart from algebras of bounded operators a
family $\cal Y$ may be not extendable to an algebra, since a
distributive law generally may be invalid, as it was seen above.
Another property is the following:
\par $(4)$ if $\{ \mbox{}_bT: ~ b\in \Psi \} $ is a net of bounded
operators in $L_q(X)$ so that $\lim_b \mbox{ }_bT=T$ in the strong
operator topology and $\mbox{}_bTA\subset B\mbox{ }_bT$ for each $b$
in a directed set $\Psi $, where $B$ is a closed operator, then
$TA\subseteq BT$.
\par Indeed, if $x\in {\cal D}(A)$, then $\mbox{}_bTx\in {\cal
D}(B)$  and $\lim_b B\mbox{ }_bTx=\lim_b \mbox{ }_bTAx=TAx$ and
hence $\lim_b \mbox{ }_bTx=Tx$. Then $Tx\in {\cal D}(B)$ and
$BTx=TAx$, since the operator $B$ is closed, from which property
$(4)$ follows. We sum up these properties as the lemma.
\par {\bf 19. Lemma.} {\it Suppose $A$ is a closed operator acting in a
Hilbert space $X$ over the Cayley-Dickson algebra ${\cal A}_v$ and
$CA\subseteq  AC$ for every operator $C$ in a self-adjoint subfamily
$\cal Y$ of $L_q(X)$, where $2\le v$. Then $TA\subseteq AT$ for each
operator $T$ in the von Neumann algebra over ${\cal A}_v$ generated
by $\cal Y$.}
\par {\bf 20. Definition.} When $A$ is a closed $\bf R$ homogeneous
${\cal A}_v$ additive operator in a Hilbert space $X$ over the
Cayley-Dickson algebra ${\cal A}_v$, $~2\le v$, and $E$ is an ${\cal
A}_v$ graded projection on $X$ so that $EA\subseteq  AE$ and $AE$ is
a bounded operator on $X$, we say that $E$ is a bounding ${\cal
A}_v$ graded projection for $A$.
\par An increasing sequence $\{ \mbox{}_nE: ~ n\in {\bf N} \} $ such
that each $\mbox{}_nE$ is a bounding ${\cal A}_v$ graded projection
for $A$ and $\vee_n \mbox{ }_nE=I$ will be called a bounding ${\cal
A}_v$ graded sequence for $A$.
\par {\bf 21. Lemma.} {\it Let $E$ be a bounding ${\cal A}_v$ graded
projection for a closed densely defined operator $T$ in a Hilbert
space over the Cayley-Dickson algebra ${\cal A}_v$, $~2\le v$. Then
$E$ is bounding for $T^*$, $T^*T$ and $TT^*$ and
\par $(1)$ $(TE)^* = T^*
E^*.$ \\ Moreover, if $\{ \mbox{}_nE: ~ n \} $ is a bounding
sequence for $T$, then $\bigcup_n \mbox{ }_nE(X)$ is a core for $T$
and $T^*$ and $T^*T$ and $TT^*$.}
\par {\bf Proof.} The conditions of this lemma mean that
$ET\subseteq  TE$ and $TE$ is bounded, hence the operator $ET$ is
pre-closed and densely defined and bounded.
\par Therefore, the operator $ET$ has the closure
$ TE$, also
\par $(2)$ $(TE)^* = (ET)^*$ is
closed and hence the operators $(TE)^*$ and $(ET)^*$ are closed by
Theorem I.3.34 \cite{ludopalglamb}. For each vectors $x\in E(X)$ and
$y\in {\cal D}(T)$ the equality $<Ty;x> = <y;(E^*T)^*x>$ is
satisfied so that $x\in {\cal D}(T^*)$ and $T^*x=(ET)^*x$,
consequently, $T^*E=(E^*T)^*E$. At the same time we have
$(I-E){\overline {E^*T}}=0$ and hence $(E^*T)^* = {\overline
{(E^*T)}}^*=(E^*T)^*E=T^*E$. Therefore, $E^*T^*\subseteq (TE)^*$ and
$E$ is bounding for $T^*$ in accordance with Equality $(2)$.
Analogously $E$ is bounding for $T^*T$ and similarly for $TT^*$,
since $E(T^*T)\subseteq (T^*T)E$.
\par This implies that $\{ \mbox{}_nE: ~ n \} $ is a bounding
sequence for $T^*$, $T^*T$ and $TT^*$ if it is such for $T$. Then
$\lim_n \mbox{ }_nEx=x$, $~\mbox{}_nEx\in {\cal D}(T)$ and $\lim_n
T\mbox{ }_nEx=  \lim_n \mbox{}_nETx=Tx$ \\
for every $x\in {\cal D}(T)$. Thus $\bigcup_{n=1}^{\infty } \mbox{
}_nE({\cal D}(T))$ is a core for $T$, consequently,
$\bigcup_{n=1}^{\infty } \mbox{ }_nE(X)$ is a core for $T$ and
$T^*$, $~T^*T$ and $TT^*$ as well, since $\mbox{}_nE(X)\subseteq
{\cal D}(T)$ for each natural number $n$.
\par {\bf 22. Remark.} Mention that in accordance with Theorem 2.28
\cite{lujmsalop} an ${\cal A}_v$ graded projection operator can be
chosen self-adjoint $E^*=E$ for a quasi-commutative von Neumann
algebra over ${\cal A}_v$, since $E$ corresponds to a characteristic
function $e$ which is real-valued.
\par {\bf 23. Theorem.} {\it  Suppose that $\sf A$ is a quasi-commutative
von Neumann algebra over either the quaternion skew field or the
octonion algebra ${\cal A}_v$, $~2\le v\le 3$, acting on a Hilbert
space $X$ over ${\cal A}_v$, also $T$ and $B\eta \sf A$. Then
\par $(1)$ each finite set of operators affiliated with $\sf A$ has
a common bounding sequence in $\sf A$; \par $(2)$ an operator $B+T$
is densely defined and pre-closed and its closure is $B{\hat +}T\eta
\sf A$; \par $(3)$ $BT$ is densely defined and pre-closed with the
closure $B{\hat \cdot }T\eta \sf A$ affiliated with $\sf A$; \par
$(4)$ $\mbox{}^jB{\hat \cdot }\mbox{ }^kT = (-1)^{\kappa
(j,k)}\mbox{ }^kT{\hat \cdot }\mbox{ }^j B$ for each $j, k$, also
$B^* B = B^* {\hat \cdot } B = B B^*$;
\par $(5)$ $((bI)B{\hat +} T)^* = B^*(b^*I){\hat +} T^* $ for each
quaternion or octonion number $b\in {\cal A}_v$; \par $(6)$ $(B{\hat
\cdot }T)^*=T^*{\hat \cdot }B^*$; \par $(7)$ if $B\subseteq T$, then
$B=T$; if $B$ is symmetric, then $B^*=B$;
\par $(8)$ the family ${\cal N}({\sf A})$ of all operators affiliated with
$\sf A$ forms a quasi-commutative $*$-algebra with unit $I$ under
the operations of addition $\hat +$ and multiplication $\hat \cdot $
given by $(2,3)$.}
\par {\bf Proof.} Take an arbitrary unitary operator $U$ in the
super-commutant ${\sf A}^{\star }$ of ${\sf A}$ (see \S II.2.71
\cite{ludopalglamb}). From $TU=UT$ it follows that $T^*U=UT^*$,
hence $T^*\eta \sf A$. Then $(T^*T)U= U(T^*T)$, since $U^*U=UU^*=I$,
consequently, $(T^*T)\eta {\sf A}$. For an ${\cal A}_v$ graded
projection $E$ in $\sf A$ an operator $(2E-I)$ is unitary such that
\par
$T(2E-I)=(2E-I)T$, \\
i.e. $T(2E-I)x=(2E-I)Tx$ for each $x\in {\cal D}(T)\subset X$, hence
$ET\subseteq TE$. In view of Theorem I.3.34 \cite{ludopalglamb} the
operator $T^*T$ is self-adjoint. Take an ${\cal A}_v$ graded
resolution $\{ \mbox{}_bE: ~ b \} $ of the identity for $T^*T$ and
put $\mbox{}_nF = \mbox{}_nE- \mbox{ }_{-n}E$ for each natural
number $n$. From Theorem 16 the inclusion $\mbox{}_bE\in \sf A$
follows. The operator $T^*T\mbox{ }_nF$ is bounded and everywhere
defined, consequently, the operator $T\mbox{ }_nF$ is everywhere
defined and closed, since $T$ is closed and $\mbox{}_nF$ is bounded.
In accordance with the closed graph theorem 1.8.6 \cite{kadring} for
$\bf R$ linear operators one gets that $T\mbox{ }_nF$ is bounded. It
can lightly be seen also from the estimate \par $ \| T\mbox{ }_nFx
\| ^2 = <\mbox{}_nFx; T^*T\mbox{ }_nFx>\le \| T^*T\mbox{ }_nF \| \|
x \| ^2$. A sequence $ \{ \mbox{}_nF: ~ n \} $ of projections is
increasing with least upper bound $I$ for which
$\mbox{}_nFT\subseteq T\mbox{}_nF$. Therefore, the limits exist:
\par $\lim_n \mbox{ }_nFx=x$ and $\lim_n T\mbox{}_nFx
= \lim_n \mbox{}_nFT x = Tx$ \\ for each vector $x$ in the domain
${\cal D}(T)$, consequently, $\bigcup_{n=1}^{\infty } \mbox{
}_nF(X)$ is a core for $T$ and $ \mbox{}_nF$ is a bounding ${\cal
A}_v$ graded sequence in $\sf A$ for the operator $T$.
\par Let now $\{ \mbox{}_nE: ~ n \} $ be a bounding sequence in $\sf
A$ for $ \{ \mbox{}_pT: ~ p=1,...,m-1 \} \subset \sf A$ and let
$\mbox{}_nF$ be a bounding sequence in $\sf A$ for $\mbox{}_mT\in
\sf A$. Then $ \{ \mbox{}_nE\mbox{ }_nF: ~ n \} $ is a bounding
sequence in $\sf A$ for $\mbox{}_1T,...,\mbox{}_mT$, particularly,
$\bigcup_{n=1}^{\infty } \mbox{}_nE\mbox{ }_nF(X)$ is a common core
for $\mbox{}_1T,...,\mbox{}_mT$. This implies that two operators
$T+B$ and $T^*+B^*$ are densely defined, but $T^*+B^*\subseteq
(T+B)^*$, consequently, $(T+B)^*$ is densely defined and $T+B$ is
pre-closed (see also Theorem I.3.34 \cite{ludopalglamb}).
\par As soon as $ \{ \mbox{}_nE: ~n \} $ is a bounding sequence in
$\sf A$ for $T$, $B$, $T^*$ and $B^*$, the inclusions are satisfied
$\mbox{}_nET\subseteq  T\mbox{}_nE$ and $\mbox{}_nEB\subseteq
 B\mbox{}_nE$ and hence
$\mbox{}_nE(TB)\subseteq (TB)\mbox{}_nE$ and $T\mbox{}_nEB\mbox{
}_nE\subseteq TB\mbox{}_nE\mbox{ }_nE$. The operators $T\mbox{ }_nE$
and $B\mbox{ }_nE$ are bounded and defined everywhere, consequently,
\par $(9)$ $T\mbox{}_nEB\mbox{ }_nE=
TB\mbox{}_nE$ \\  (see also Theorem 2.28 \cite{lujmsalop}). This
means that $ \{ \mbox{}_nE: ~ n \} $ is a bounding sequence for $TB$
and analogously for $BT$ and $B^*T^*$. That is, the operator
$B^*T^*$ is densely defined. Since $B^*T^*\subseteq (TB)^*$, one
gets that the operator $(TB)^*$ is densely defined and $TB$ is
pre-closed. From Formula $(9)$ we infer that
$$\mbox{}^kT\mbox{ }^jB\mbox{ }_nE = (-1)^{\kappa (j,k)}\mbox{ }^jB
\mbox{ }^kT \mbox{ }_nE$$ for each $j, k$, since $\sf A$ is
quasi-commutative. Therefore, the operators $T{\hat \cdot }B$ and
$B{\hat \cdot }T$ agree on their common core $\bigcup_{n=1}^{\infty
} \mbox{ }_nE(X)$ and inevitably $$(10)\quad \mbox{}^jB{\hat \cdot
}\mbox{ }^kT =(-1)^{\kappa (j,k)} \mbox{ }^kT{\hat \cdot } \mbox{
}^jB.$$ The operators $T^*T$ and $TT^*$ are self-adjoint,
consequently, $T^*T=T^*{\hat \cdot }T$ and $TT^*=T{\hat \cdot }T^*$
and $T^*{\hat \cdot }T=  T{\hat \cdot }T^*$. Then $U^*x$ and $Ux\in
{\cal D}(T+B)$ for each $x\in {\cal D}(T)\cap {\cal D}(B)={\cal
D}(T+B)$, consequently, $U{\cal D}(T+B)={\cal D}(T+B)$ and $(T+B)U=
U (T+B)$. Thus $(T{\hat +} B)U = U (T{\hat +} B)$, as well as
$(T{\hat +}B)\eta \sf A$.
\par For each vector $y\in {\cal D}(TB)$ the inclusions follow $y\in
{\cal D}(B)$ and $By\in {\cal D}(T)$. Therefore, $Uy\in {\cal D}(B)$
and $BUy = UB y\in {\cal D}(T)$, consequently, $Uy\in {\cal D}(TB)$.
Then we get $U{\cal D}(TB)={\cal D}(TB)$, since $U^*y\in {\cal
D}(TB)$. But $(TB)Uy= U(TB)y$, hence $(T{\hat \cdot }B)Uy= U(T{\hat
\cdot }B)y$, and inevitably $(T{\hat \cdot }B)\eta \sf A$. \par
Having a bounding ${\cal A}_v$ graded sequence $ \{ \mbox{}_nE: ~ n
\} $ for $T$ and $T^*$ one gets $\mbox{}_nET^*\subseteq  T^*\mbox{
}_nE$ and $\mbox{}_nE^*T^*\subseteq (T\mbox{ }_nE)^*$, consequently,
$T^*\mbox{ }_nE$ and $(T\mbox{ }_nE)^*$ are bounded everywhere
defined extensions of operators $ \mbox{}_nET^*$ and $\mbox{
}_nE^*T^*$ respectively. In view of Lemma 21 the equality $(T\mbox{
}_nE)^* = \mbox{}_nE^*T^*$ follows. \par We now consider a bounding
${\cal A}_v$ graded sequence $ \{ \mbox{}_nE: ~ n \} $ for $T$, $~
T^*$, $~B$, $~B^*$, $ ~ ((bI)T{\hat +}B)$, $ ~ ((bI)T{\hat +}B)^*$,
$ ~ (T{\hat \cdot }B)$, $ ~ (T{\hat \cdot }B)^*$ and $T^*{\hat \cdot
}B^*$. From the preceding demonstration we get the equalities
\par $ ~ (T^* (b^*I) {\hat +}B^*)\mbox{ }_nE =
(T^* (b^*I)) \mbox{ }_nE {\hat +} B^*\mbox{ }_nE =  [((bI)T{\hat
+}B)\mbox{ }_nE^*]^*= ((bI)T{\hat +}B)^*\mbox{ }_nE$ and
\par $(T{\hat \cdot }B)^*\mbox{ }_nE = [(T{\hat \cdot }B)
\mbox{ }_nE^*]^*=
(B^* {\hat \cdot }T^*)\mbox{ }_nE$ \\
due to $(9,10)$. The operators $((bI)T{\hat +}B)^*$ and
$(T^*(b^*I){\hat +}B^*)$ agree on their common core
$\bigcup_{n=1}^{\infty } \mbox{ }_nE(X)$, consequently, $((bI)T{\hat
+}B)^* = (T^*(b^*I){\hat +}B^*)$. Analogously we infer that $(T{\hat
\cdot }B)^* = B^*{\hat \cdot }T^*$.
\par Then $T\mbox{ }_nE\subseteq B\mbox{ }_nE$ and hence
$T\mbox{ }_nE = B\mbox{ }_nE$ as soon as $T\subseteq B$ and $ \{
\mbox{}_nE: ~ n \} $ is a bounding ${\cal A}_v$ graded sequence in
$\sf A$ for $T$ and $B$. Therefore, the operators $T$ and $B$ are
consistent on their common core $\bigcup_{n=1}^{\infty } \mbox{
}_nE(X)$ and hence $T=B$. If $T$ is symmetric, then $T\subseteq T^*$
and from the preceding conclusion one obtains $T=T^*$.
\par For any three operators $T, B, C \in {\cal N}({\sf A})$
we take a common bounding ${\cal A}_v$ graded sequence $ \{
\mbox{}_nE: ~ n \} $ and get
\par $(T{\hat \cdot }B){\hat \cdot }C =
 T{\hat \cdot }(B{\hat \cdot }C)$, since
\par $[(T{\hat \cdot }B){\hat \cdot }C]\mbox{ }_nE=
 [T{\hat \cdot }(B{\hat \cdot
}C)]\mbox{ }_nE$ \\ for all $n$ (see also \cite{schafb}). From this
Statement $(8)$ of the theorem follows.
\par {\bf 24. Lemma.} {\it Suppose that $ \{ \mbox{}_nF: ~ n \} $
is a bounding ${\cal A}_v$ graded
sequence for the closed operator $T$ on a Hilbert space $X$ over the
Cayley-Dickson algebra ${\cal A}_v$, $~2\le v$, and $T\mbox{ }_nF$
is normal for each natural number $n$. Then $T$ is normal.}
\par {\bf Proof.} In view of Lemma 21 we have
the equalities $(T\mbox{ }_nF)^* = \mbox{}_nF^*T^*$ and
$\mbox{}_nFT^*=  T^* \mbox{}_nF$ and $T^*T\mbox{ }_nF =T^*\mbox{
}_nFT\mbox{ }_nF$ and $(T^* T)\mbox{ }_nF = (TT^*)\mbox{ }_nF$.
Therefore, the self-adjoint operators $T^*T$ and $TT^*$ agree on
their common core $\bigcup_{n=1}^{\infty } \mbox{ }_nF(X)$,
consequently, $T^*T = TT^*$, i.e. the operator $T$ is normal.
\par {\bf 25. Remark.}  The condition $T^*T=TT^*$
is equivalent to
\par $\sum_{j,k} [(T^*)^{i_k}T^{i_j} - T^{i_j}(T^*)^{i_k}]=0$
on ${\cal D}(T^*T)={\cal D}(TT^*)$, since $\sum_j \pi ^j=I$ (see \S
2). \par If $BT\subseteq TB$ in a Hilbert space $X$ over the
Cayley-Dickson algebra, then ${\bf B}{\bf T}\subseteq {\bf T}{\bf
B}$ for ${\bf B}={\bf i} B={{0 ~~ B}\choose {-B ~ 0}}$ and ${\bf
T}={\bf i}T={{0 ~~ T}\choose {-T ~ 0}}$ defined on ${\cal D}({\bf
B})={\cal D}(B)\oplus {\cal D}(B){\bf i}$ and ${\cal D}({\bf T}) =
{\cal D}(T)\oplus {\cal D}(T){\bf i}$ with ${\bf i} ={{0 ~~
1}\choose {-1 ~ 0}}$ (see also \S 18).
\par {\bf 26. Lemma.} {\it Let $BT\subseteq  TB$
and ${\cal D}(T)\subseteq {\cal D}(B)$, let also $T$ be a
self-adjoint operator and let $B$ be a closed operator in a Hilbert
space $X$ over either the quaternion skew field or the octonion
algebra ${\cal A}_v$, $~2\le v\le 3$. Then $\mbox{}_bEB\subseteq
B\mbox{ }_bE$ for each $\mbox{}_bE$ in the spectral ${\cal A}_v$
graded resolution $ \{ \mbox{}_bE: ~ b \} $ of $T$.}
\par {\bf Proof.} There is the decomposition $alg_{{\cal A}_v}
(I,B,T,{\bf i}I, {\bf i}T, {\bf i}T) = alg_{{\cal A}_v} (I,B,T)
\oplus alg_{{\cal A}_v} (I,B,T){\bf i}$ and this family is defined
in the Hilbert space $X\oplus X{\bf i}$. From Formula 18$(3)$ we
have the inclusion $(BT+B{\bf i}I)\subseteq B(T+{\bf i}I)$. Now we
consider these operators with domains in $X\oplus X{\bf i}$ denoted
by ${\cal D}(T)$, ${\cal D}(B)$, etc. Take an arbitrary vector $x\in
{\cal D}(B(T+{\bf i}I))$, then $x\in {\cal D}(T)$ and $Tx+{\bf
i}x\in {\cal D}(B)$. But from the suppositions of this lemma we have
the inclusion ${\cal D}(T)\subseteq {\cal D}(B)$, consequently,
$Tx\in {\cal D}(B)$ and hence $x\in {\cal D}(BT+B({\bf i}I))$ and
$B(T+{\bf i}I)x=BTx+B({\bf i}x)$. Thus $B(T+{\bf i}I)\subseteq
BT+B({\bf i}I)$, consequently, $B(T+{\bf i}I)=BT+B({\bf i}I)$.
\par Denote by $Q_-$ and $Q_+$ the bounded everywhere defined
inverses to $(T-{\bf i}I)$ and $(T+{\bf i}I)$ respectively. In view
of 18$(1-3)$ and the preceding proof we infer \par $Q_{\pm }B=Q_{\pm
}B(T\pm {\bf i}I)Q_{\pm }=Q_{\pm }(BT\pm B({\bf i}I))Q_{\pm
}\subseteq Q_{\pm }(TB\pm B({\bf i}I))Q_{\pm }=Q_{\pm }(T\pm {\bf
i}I)BQ_{\pm }$ and hence $(Q_{\pm })B\subseteq B(Q_{\pm })$. From \S
10 one has $Q_+ = (Q_-)^*$ and applying Lemma 19 one gets
$QB\subseteq BQ$ for each element $Q$ in the von Neumann algebra
over ${\cal A}_v$ generated by $Q_+$ and $Q_-$. To an ${\cal A}_v$
graded projection operator $\mbox{}_bE$ on $X$ an operator
$\mbox{}_bE\oplus \mbox{}_bE$ on $X\oplus X{\bf i}$ corresponds.
Particularly, this means the inclusion $\mbox{}_bEB\subseteq B\mbox{
}_bE$ in $X$ for each $b$.
\par {\bf 27. Theorem.} {\it Suppose that $T$ is an operator on a
Hilbert space $X$ over either the quaternion skew field or the
octonion algebra ${\cal A}_v$, $~2\le v\le 3$. An operator $T$ is
normal if and only if it is affiliated with a quasi-commutative von
Neumann algebra $\sf A$ over ${\cal A}_v$. Moreover, there exists
the smallest such algebra $\mbox{}_0{\sf A}$.}
\par {\bf Proof.} In view of Theorem 23 if an operator is
affiliated with a quasi-commutative von Neumann algebra $\sf A$ over
either the quaternion skew field or the octonion algebra ${\cal
A}_v$ with $2\le v \le 3$, then it is normal. Suppose that an
operator $T$ is normal. Then Lemma 26 is applicable, since
$TT^*T=T^*TT$ and ${\cal D}(T^*T)\subseteq {\cal D}(T)$. \par For an
${\cal A}_v$ graded spectral resolution $\{ \mbox{}_bE: ~ b\in {\bf
R} \} $ this implies that $\mbox{}_bET\subseteq T\mbox{ }_bE$ for
each $b$, hence $\mbox{}_nFT\subseteq T\mbox{ }_nF$ for each natural
number $n$, where $\mbox{}_nF = \mbox{}_nE - \mbox{ }_{-n}E$. Then
we also have $T^*T^*T = T^*TT^*$ and ${\cal D}(T^*T)={\cal
D}(TT^*)\subseteq {\cal D}(T^*)$, consequently,
$\mbox{}_nFT^*\subseteq T^*\mbox{ }_nF$. From \S 23 it follows that
the operators $T\mbox{ }_nF$ and $\mbox{}_nFT$ are bounded, since
such the operator $T^*T\mbox{ }_nF = TT^*\mbox{ }_nF$ is. But then
$\mbox{}_nF^*T^*\subseteq (T\mbox{}_nF)^*$ and hence both operators
$(T\mbox{ }_nF)^*$ and $T^*\mbox{ }_nF$ are bounded extensions of
the densely defined operator $\mbox{ }_nFT^*$ when choosing
$\mbox{}_bE$ and hence $\mbox{}_nF$ self-adjoint for each $b$ and
$n$. Thus one gets the equalities $(T\mbox{ }_nF)^* = T^*\mbox{
}_nF$ and $(T^*\mbox{ }_nF)^* = T\mbox{ }_nF$. On the other hand,
there are the inclusions $T\mbox{ }_nFT\mbox{ }_mF \subseteq
TT\mbox{ }_nF$ and $T\mbox{ }_mFT\mbox{ }_nF\subseteq TT\mbox{ }_nF$
for each $n\le m$. \par The operators $T\mbox{ }_nFT\mbox{ }_mF$ and
$T\mbox{ }_mFT\mbox{ }_nF$ are everywhere defined, consequently,
$T\mbox{ }_nFT\mbox{ }_mF = TT\mbox{ }_nF= T\mbox{ }_mFT\mbox{
}_nF.$ There are the equalities $T^* \mbox{ }_mF T \mbox{ }_nF
=T^*T\mbox{ }_nF=TT^*\mbox{ }_nF=T\mbox{ }_nFT^*\mbox{ }_mF$, hence
$$\mbox{}_0{\sf A} := cl [ alg_{{\cal A}_v} \{ \mbox{}_nF, ~ T\mbox{
}_nF, ~ T^*\mbox{ }_nF: ~ n\in {\bf N} \} ]$$ is a quasi-commutative
von Neumann algebra over the algebra ${\cal A}_v$, since there is
the inclusion $alg_{{\cal A}_v} \{ \mbox{}_nF, ~ T\mbox{ }_nF, ~
T^*\mbox{ }_nF: ~ n\in {\bf N} \} \subset L_q(X)$.
\par Moreover, one has that $\bigcup_{n=1}^{\infty } \mbox{ }_nF(X)=:
\cal Y$ is a core for $T$, since $\vee_{n=1}^{\infty } \mbox{
}_nF=I$ and $\mbox{}_nFT\subseteq T\mbox{ }_nF$. If $U$ is a unitary
operator in $(\mbox{}_0{\sf A})^{\star }$ and $x\in \cal Y$, then
the equalities $TUx=TU\mbox{ }_nFx=T\mbox{ }_nFUx = UT\mbox{ }_nFx =
UTx$ are fulfilled for some natural number $n$. In accordance with
Remark 2 $T\eta \mbox{ }_0{\sf A}$ and  $T^*\eta \mbox{ }_0{\sf A}$
also. If $T\eta {\cal R}$, then $T^*\eta {\cal R}$ and $T^*T\eta
{\cal R}$ are affiliations  as well. From Note 15 it follows that
the self-adjoint operator $T^*T$ generates a quasi-commutative von
Neumann algebra $\sf A$ contained in $\cal R$, consequently,
$\mbox{}_nF\in \cal R$ and hence $T\mbox{ }_nF, ~ T^*\mbox{ }_nF\in
\cal R$. This implies the inclusion $\mbox{}_0{\sf A}\subset \cal
R$.
\par {\bf 28. Definition.} The algebra $\mbox{}_0\sf A$ from the
preceding section will be called the von Neumann algebra over the
Cayley-Dickson algebra ${\cal A}_v$ with $2\le v$ generated by the
normal operator $T$.
\par {\bf 29. Theorem.} {\it Suppose that $\sf A$ is a quasi-commutative
von Neumann algebra over either the quaternion skew field or the
octonion algebra ${\cal A}_v$, $~2\le v\le 3$, also $\phi $ is an
isomorphism of $\sf A$ onto $C(\Lambda ,{\cal A}_v)$, where $\Lambda
$ is a compact Hausdorff topological space, $T\eta \sf A$. Then
there exists a unique normal function $\phi (T)$ on $\Lambda $ so
that $\phi (TE)=\phi (T){\hat \cdot }\phi (E)$, when $E$ is an
${\cal A}_v$ graded bounding projection in $\sf A$ for $T$, where
$(\phi (T){\hat \cdot } \phi (E))(z)=\phi (T)(z) \phi (E)(z)$, if
$\phi (T)(z)$ is defined and zero in the contrary case, $z\in
\Lambda $. If ${\cal N}(\Lambda ,{\cal A}_v)$ is the family of all
${\cal A}_v$ valued normal function on $\Lambda $ and $f, g \in
{\cal N}(\Lambda ,{\cal A}_v)$, then there are unique normal
functions ${\tilde f}$, $~sf$, $~fs$, $~f{\hat +}g$ and $f{\hat
\cdot }g$ so that ${\tilde f}(z)=\widetilde{f(z)}$,
$~(sf)(z)=sf(z)$, $~(fs)(z)=f(z)s$, $~(f{\hat +}g) (z) = f(z)+g(z)$
and $(f{\hat \cdot }g)(z) = f(z){\hat \cdot }g(z)$, when $f$ and $g$
are defined at $z\in \Lambda $, $~s\in {\cal A}_v$. Endowed with the
operations $f\mapsto {\tilde f}$ and $(s,f)\mapsto sf$ and
$(f,s)\mapsto fs$ and $(f,g)\mapsto f{\hat +}g$ and $(f,g)\mapsto
f{\hat \cdot }g$, the family ${\cal N}(\Lambda ,{\cal A}_v)$ is a
quasi-commutative algebra with unit $1$ and involution $f\mapsto
\tilde f$, it is associative over the quaternion skew field ${\bf
H}={\cal A}_2$ and alternative over the octonion algebra ${\bf
O}={\cal A}_3$. The natural extension of $\phi $ is a
$*$-isomorphism of ${\cal N}({\sf A})$ onto ${\cal N}(\Lambda ,{\cal
A}_v)$.}
\par {\bf Proof.} In view of Theorem 23 an operator $T$ affiliated
with $\sf A$ has an ${\cal A}_v$ graded bounding sequence $ \{
\mbox{}_nE: ~ n \} $ in $\sf A$. If $\phi $ has the properties
described above, if also $f$ and $g$ are normal functions defined at
a point $z$ and corresponding to $\phi (T)$, then
\par $f(z)\phi (\mbox{}_nE)(z) = \phi (T\mbox{ }_nE)(z) =g(z) \phi
(\mbox{}_nE)(z)$ \\ for every natural number $n$. Thus if $\phi
(\mbox{}_nE)(z)=1$ for some natural number $n$, then certainly
$f(z)=g(z)$. Applying Lemma 9 we obtain, that $f=g$, since the
sequence $\mbox{}_nE$ is monotone increasing to $I$, while $f$ and
$g$ agree on a dense subset of $\Lambda $. Thus an isomorphism $\phi
$ is unique.
\par On the other hand, Theorem 16 provides $\phi (T)$ with the
required properties for each $T\in {\cal Q}({\sf A})$, where ${\cal
Q}({\sf A})$ denotes a family of all self-adjoint operators
affiliated with $\sf A$. This means that $\phi $ is defined on
${\cal Q}({\sf A})$. \par For any pair of operators $T$, $B$ in
${\cal Q}({\sf A})$ it is possible to choose a bounding sequence $\{
\mbox{}_nE: ~ n \} $ for both $T$ and $B$ by Theorem 23. Therefore,
the operators $T\mbox{ }_nE$, $ ~ B\mbox{ }_nE$, $~(T+B)\mbox{ }_nE$
and $TB\mbox{ }_nE$ belong to the quasi-commutative von Neumann
algebra $\sf A$ over the Cayley-Dickson algebra ${\cal A}_v$ so that
\par $(T{\hat +}B)\mbox{ }_nE= T\mbox{ }_nE + B\mbox{ }_nE$ and
\par $(T{\hat \cdot }B)\mbox{ }_nE = TB\mbox{ }_nE.$ Therefore, we
deduce the equalities
\par $\phi (T{\hat +}B)(z)\phi (\mbox{}_nE)(z) =
\phi ((T{\hat +}B)\mbox{}_nE)(z) = \phi (T)(z) \phi (\mbox{}_nE)(z)
+ \phi (B)(z) \phi (\mbox{}_nE)(z)$ and
\par $\phi (T{\hat \cdot }B)(z)\phi (\mbox{}_nE)(z) =
\phi ((T{\hat \cdot }B)\mbox{}_nE)(z) = \phi (T)(z)\phi (B)(z) \phi
(\mbox{}_nE)(z)$, \\
when $\phi (T)$, $~\phi (B)$, $~ \phi (T{\hat +}B)$ and $\phi
(T{\hat \cdot } B)$ are defined at $z$. The pairs $[\phi (T{\hat
+}B); \phi (T)+\phi (B)]$ and $[\phi (T{\hat \cdot } B); \phi
(T)\phi (B)]$ agree on dense subsets of $\Lambda $. Therefore,
images $\phi (T{\hat +}B)$ and $\phi (T{\hat \cdot }B)$ are finite,
when $\phi (T)$ and $\phi (B)$ are defined, hence $\phi (T{\hat
+}B)$ and $\phi (T{\hat \cdot }B)$ are normal extensions of $\phi
(T)+\phi (B)$ and $\phi (T)\phi (B)$ correspondingly. Then $\phi
(T){\hat +} \phi (B)$ and $\phi (T){\hat \cdot } \phi (B)$ are
defined as normal extensions of $\phi (T)+\phi (B)$ and $\phi
(T)\phi (B)$ respectively. Each function $w\in {\cal Q}(\Lambda )$
corresponds to some operator $T\in {\cal Q}({\sf A})$ and by Theorem
16 the operations ${\hat +}$ and $\hat \cdot $ are applicable to all
functions in ${\cal Q}(\Lambda )$. An iterated application of Lemma
9 shows that ${\cal Q}(\Lambda )$ endowed with these operations is a
commutative algebra over the real field $\bf R$ with unit $1$, since
each function in ${\cal Q}(\Lambda )$ is real valued. Moreover,
$\phi $ is an isomorphism of ${\cal Q}({\sf A})$ onto ${\cal
Q}(\Lambda )$.
\par If $v$ is finite, for each $j=0,...,2^v-1$ the $\bf R$-linear
projection operator $\hat \pi ^j := \hat \pi ^j_v: {\sf A}\to {\sf
A}_ji_j$ is expressible as a sum of products with generators and
real constants due to Formulas $(1,2)$ below so that $\hat \pi
^j(A)=i_jA_j=A_ji_j$:
\par $(1)$ $\hat \pi ^j(A) = (- i_j (Ai_j) - (2^v-2)^{-1} \{ - A
+\sum_{k=1}^{2^v-1}i_k(Ai_k^*) \} )/2$
\\ for each $j=1,2,...,2^v-1$, $$(2)\quad \hat \pi ^0(A) = (A+
(2^v-2)^{-1} \{ -A + \sum_{k=1}^{2^v-1}i_k(Ai_k^*) \} )/2,$$  where
$2\le v\in \bf N$, \par $(3)$ $ A=\sum_j A_ji_j$, \\ $ ~ A_j\in {\sf
A}_j$ for each $j$, $~A\in \sf A$.
\par If $\sf A$ is embedded into $L_q(X)$ for a Hilbert space $X$
over the algebra ${\cal A}_v$, we can consider projections $\pi ^j :
{\sf A}\to {\sf A}_ji_j$ and get decomposition $(3)$ so that $A_j\in
\sf A$ for each $j$. The function $\sum_j \phi (A_j)i_j$ is defined
on $\Lambda \setminus \bigcup_j W_j$ when $\phi (A_j)$ is defined on
$X\setminus W_j$. But $\Lambda \setminus \bigcup_j W_j$ is
everywhere dense in $\Lambda $, since $|z|=\sqrt{\sum_j z_j^2}$ is
the norm on ${\cal A}_v$. If $p\in W_j$, then $$\lim_{x\in \Lambda
\setminus \bigcup_j W_j, ~ x\to p } |\phi (A)(x)| =\sqrt{\sum_j \phi
(A_j)^2(x)}=\infty .$$ That is $\sum_j \phi (A_j)i_j$ is an element
${\hat {\sum }}_j \phi (A_j) i_j= \phi (A)\in {\cal N}(\Lambda
,{\cal A}_v)$.
\par In accordance with Theorem 6.2.26 \cite{eng} a topological
space $\Lambda $ is extremely disconnected if and only if for each
pair of nonintersecting open subsets $U$ and $V$ in $\Lambda $ their
closures always $cl (U)\cap cl (V)=\emptyset $ do not intersect. By
Theorem 8.3.10 \cite{eng} if $(S,{\cal U})$ is a uniform space and
$(Y,{\cal V})$ is a complete uniform space, then each uniformly
continuous mapping $f: (P,{\cal U}_P)\to (Y,{\cal V})$, where $P$ is
an everywhere dense subset in $S$ relative to a topology induced by
a uniformity $\cal U$, has a uniformly continuous extension $f:
(S,{\cal U})\to (Y,{\cal V})$.
\par A set $\xi ^{-1}(0)$ is closed in $\Lambda $, when $\xi \in
{\cal N}(\Lambda ,{\cal A}_v)$. Indeed, a function $\xi $ is
continuous on $\Lambda \setminus W_{\xi }$, consequently, $\xi
^{-1}(0)$ is closed in $\Lambda \setminus W_{\xi }$. At the same
time the limit $$\lim_{x\to p, ~ x\in \Lambda \setminus W_{\xi }}
|\xi (x)| = \infty $$ is infinite for each $p\in W_{\xi }$, hence
$p\notin cl [\xi ^{-1}(0)]$. This fact implies that the interior $K
:= Int [\xi ^{-1}(0)]$ is clopen in $\Lambda $. Thus the function
$g$ defined such that $g(x)=1$ on $K$, while $g(x) := \xi (x)/|\xi
(x)|$ on $\Lambda \setminus (\xi ^{-1}(0)\cup W_{\xi })$, is
continuous on $[\Lambda \setminus (\xi ^{-1}(0)\cup W_{\xi })]\cup
K$. This function $g$ has a continuous extension $q\in C(\Lambda
,{\cal A}_v)$ by Theorems 6.2.26 and 8.3.10 \cite{eng}, since
$\Lambda \setminus (\{ \xi ^{-1}(0)\setminus K \} \cup W_{\xi })$ is
a dense open subset in $\Lambda $.
\par Indeed, Cantor's cube $D^{\sf m}$ is universal for all
zero-dimensional spaces of topological weight ${\sf m}\ge \aleph
_0$, where $D$ is a discrete two-element space (see Theorem 6.2.16
\cite{eng}). If ${\sf m}<\aleph _0$, then $\Lambda $ is discrete and
this case is trivial. Therefore, if ${\sf m}\ge \aleph _0$, the
compact topological space $\Lambda $ has an embedding into $D^{\sf
m}$ as a closed subset. Then a uniformity $\cal U$ compatible with
its topology can be chosen non-archimedean. That is, each
pseudo-metric $\rho $ in a family $\bf P$ of all pseudo-metrics
inducing a uniformity $\cal U$ on $\Lambda $ satisfies the
inequality \par $\rho (x,y) \le \max (\rho (x,z); ~ \rho (z,y))$ for
each $x, y, z\in \Lambda $.
\par The function $g(x)$ has values in the unit sphere in ${\cal
A}_v$, which is closed in the Cayley-Dickson algebra ${\cal A}_v$
and hence is complete. If $z\in (\xi ^{-1}(0)\cup W_{\xi })\setminus
K$ and $\{ x_{\alpha }: \alpha \in \Sigma \} \subset [\Lambda
\setminus (\xi ^{-1}(0)\cup W_{\xi })]\cup K$ is a Cauchy net
converging to $z$, then $\{ x_{\alpha }: \alpha \in \Sigma \}$ can
be chosen such that the limit $\lim_{\alpha } g(x_{\alpha })$
exists, where $\Sigma $ is a directed set. For any other Cauchy net
$ \{ y_{\alpha } : \alpha \in \Sigma \} \subset [\Lambda \setminus
(\xi ^{-1}(0)\cup W_{\xi })]\cup K$ converging to $z$ the limit of
the net $ \{ g(y_{\alpha }) : \alpha \in \Sigma \} $ exists and will
be the same for $ \{ g(x_{\alpha }) : \alpha \in \Sigma \} $, since
each pseudo-metric $\rho \in \bf P$ is non-archimedean and the
function $g$ is continuous on $[\Lambda \setminus (\xi ^{-1}(0)\cup
W_{\xi })]\cup K$ (see also Propositions 1.6.6 and 1.6.7 and Theorem
8.3.20 \cite{eng}).
\par For a point $x\in \Lambda
\setminus W_{\xi }$ the equality $q(x)|\xi (x)|=\xi (x)$ is
fulfilled. Since $|\xi |\in {\cal Q}(\Lambda )$, one gets also
normal extensions $q_j{\hat \cdot } |\xi |$ for $\xi _j$ defined on
$\Lambda \setminus W_{\xi }$. Thus the functions $\xi , \xi _j$ are
in ${\cal N}(\Lambda ,{\cal A}_v)$ for each $j$. Choosing $A_j\in
{\cal Q}({\sf A})$ as $\phi (A_j)=\xi _j$ leads to the equalities
$\phi ({\hat \sum }_j A_ji_j) = {\hat \sum }_j \xi _ji_j=\xi $. Thus
this function $\phi $ maps ${\cal N}({\sf A})$ onto ${\cal
N}(\Lambda ,{\cal A}_v)$.
\par As soon as functions $f, g \in {\cal N}(\Lambda ,{\cal A}_v)$ are
defined on $\Lambda \setminus W_f$ and $\Lambda \setminus W_g$
respectively, then components $f_j$ and $g_j$ have normal extensions
for each $j$ as follows from the proof above. Therefore, their sum
$f+g$ and product $fg$ defined on $\Lambda \setminus (W_f\cup W_g)$
have the normal extensions ${\hat \sum }_j (f_j {\hat +} g_j) i_j =
f{\hat +}g$ and ${\hat \sum}_{j,k} (f_j {\hat {\cdot }} g_k) i_ji_k
= f{\hat {\cdot }}g$. \par The set ${\cal N}(\Lambda ,{\cal A}_v)$
supplied with these operations is the algebra over the
Cayley-Dickson algebra ${\cal A}_v$ with $2\le v$. Moreover, the
algebra ${\cal N}(\Lambda ,{\cal A}_v)$ is quasi-commutative by its
construction and with unit $1$ and adjoint operation $f\mapsto
{\tilde f}$. This algebra ${\cal N}(\Lambda ,{\cal A}_v)$ is
associative over the quaternion skew field ${\bf H}={\cal A}_2$ and
alternative over the octonion algebra ${\bf O}={\cal A}_3$, since
the multiplication of functions is defined point-wise, while the
quaternion skew field ${\bf H}$ is associative and the octonion
algebra ${\bf O}$ is alternative. Therefore, the mapping $\phi $ has
an extension up to a $*$-isomorphism of ${\cal N}({\sf A})$ onto
${\cal N}(\Lambda ,{\cal A}_v)$. For $T\eta \sf A$ and a bounding
${\cal A}_v$ graded projection $E$ in $\sf A$ for $T$ we get $\phi
(TE) = \phi (T ~ {\hat {\cdot }} ~ E) = \phi (T) ~ {\hat {\cdot }} ~
\phi (E)$.

\par {\bf 30. Definition.} The spectrum $sp (T)$ of a closed densely
defined operator $T$ on a Hilbert space $X$ over the Cayley-Dickson
algebra ${\cal A}_v$ with $2\le v$ is the set of all those
Cayley-Dickson numbers $z\in {\cal A}_v$ for which an operator
$(T-zI)$ is not a bijective $\bf R$ linear ${\cal A}_v$ additive
mapping of ${\cal D}(T)$ onto $X$, where ${\cal D}(T)$ as usually
denotes a ${\cal A}_v$ vector domain of definition of $T$.

\par {\bf 31. Remark.} For an operator $T$ from Definition 30
if $z\notin sp (T)$, then $(T-zI)$ is bijective from ${\cal D}(T)$
onto $X$ and has an $\bf R$ linear ${\cal A}_v$ additive inverse
$B=(T-zI)^{-1}: X\to {\cal D}(T)$. The graph of $B$ is closed, since
it is such for $(T - zI)$. In accordance with the closed graph
theorem 1.8.6 \cite{kadring} this inverse operator $B$ is bounded.
Thus we get $z\notin sp (T)\Leftrightarrow (T-zI)$ has a bounded
inverse from $X$ onto ${\cal D}(T)$. If $T\eta \sf A$ for some
quasi-commutative von Neumann subalgebra in $L_q(X)$ and $z\notin sp
(T)$, then $B\in \sf A$. From the boundedness of $B$ and closedness
of $(T-zI)$ it follows that \par $I=(T-zI)B= (T-zI)~ {\hat \cdot } ~
B = B ~ {\hat \cdot } ~(T-zI)$. \par Therefore, $z\notin sp (T)$ is
equivalent to: $(T-zI)$ has an inverse $B$ in the algebra ${\cal
N}(\Lambda ,{\cal A}_v)$ and $B\in \sf A$.

\par {\bf 32. Proposition.} {\it Let $T$ be a normal operator
affiliated with a quasi-commutative von Neumann algebra $\sf A$ over
either the quaternion skew field or the octonion algebra ${\cal
A}_v$ with $2\le v\le 3$. Then $sp (T)$ coincides with the range of
$\phi (T)$, where $\phi $ denotes the isomorphism of ${\cal N}({\sf
A})$ onto ${\cal N}(\Lambda ,{\cal A}_v)$ extending the isomorphism
of $\sf A$ with $C(\Lambda ,{\cal A}_v)$ (see \S 29).}
\par {\bf Proof.} In accordance with Definition 30 $z\notin sp (T)$
if and only if there is a bounded $B$ inverse to $(T-zI)$. For a
unitary operator $U$ in ${\sf A}^{\star }$ one has the equality
$U^*(T-zI)U=(T-zI)$, since $(T-zI)\eta \sf A$. Therefore, $U^*BU=B$
for each such unitary operator $U$ and $B\in \sf A$. But an operator
$B\in \sf A$ exists so that $(T-zI) \mbox{ }^{ \hat{.}} \mbox{ }
B=I$. Thus the equality $\phi (T-zI) \mbox{ }^{ \hat{.}} \mbox{ }
\phi (B) = 1$ is satisfied if and only  if $z\notin sp (T)$. This
means that there exists such $\phi (B)\in C(\Lambda ,{\cal A}_v)$ if
and only if a Cayley-Dickson number $z$ does not belong to the range
of $\phi (T)$, hence $sp (T)$ is the range of $\phi (T)$.
\par {\bf 33. Definition.} A self-adjoint operator $T$ is called
positive when \par $<Tx;x> ~ \ge 0$ \\ for each vector $x\in {\cal
D}(T)$ in its domain.
\par {\bf 34. Note.}  A function $f-z1$ for $f\in {\cal N}(\Lambda
,{\cal A}_v)$ and $z\in {\cal A}_v$ has not an inverse in ${\cal
N}(\Lambda ,{\cal A}_v)$ if and only if $f-z1$ vanishes on some
non-void clopen subset in $\Lambda $, where $2\le v \le 3$.
Considering this in ${\cal N}({\sf A})$ we get a non-zero ${\cal
A}_v$ graded projection operator $F\in \sf A$ so that $(T-zI) ~
{\hat \cdot } ~ F =0$. The latter is equivalent to the existence of
a non-zero vector $x$ on which $(T-zI)x=0$. In this situation one
says that $z$ is in the point spectrum of $T$. Thus the spectrum of
$T$ relative to ${\cal N}({\sf A})$ is its point spectrum.
\par {\bf 35. Proposition.} {\it Let $T$ be a self-adjoint operator
in a Hilbert space $X$ over either the quaternion skew field or the
octonion algebra. Then $T$ is positive if and only if $z\ge 0$ for
each $z\in sp (T)$.}
\par {\bf Proof.} In view of Lemma 10 it is sufficient to consider
the variant, when $T$ is affiliated with a quasi-commutative von
Neumann algebra $\sf A$. There exists an isomorphism $\phi $ of
${\cal N}({\sf A})$ onto ${\cal N}(\Lambda ,{\cal A}_v)$ extending
an isomorphism of $\sf A$ onto $C(\Lambda ,{\cal A}_v)$ (see Theorem
29), where $\Lambda $ is an extremely disconnected compact Hausdorff
topological space. Whenever $\phi (T)$ is defined  on $\Lambda
\setminus W_{\phi (T)}$ and $\phi (T)(y)<0$ for some point $y\in
\Lambda \setminus W_{\phi (T)}$, there exists a non-void clopen
subset $P$ containing $y$ and a negative constant $b<0$ so that
$P\subset \Lambda \setminus W_{\phi (T)}$ and $\phi (T)(x)<b<0$ for
each $x\in P$. If $F$ is an ${\cal A}_v$ graded projection in $\sf
A$ corresponding to the characteristic function $\chi _P$ of $P$,
then $F$ is a non-zero ${\cal A}_v$ graded bounding projection for
the operator $T$ so that $TF\le b F$. For each unit vector $z\in
{\cal R}(F)$ in the range of $F$, the inequality $<Tz;z> ~ \le b< 0$
would be fulfilled, hence $T$ would be not positive. This implies
the inequality $sp (T)\ge 0$ if $T\ge 0$.
\par As $\phi (T)$ has range consisting of non-negative real numbers,
its  positive square root $g$ is a normal self-adjoint mapping on
$\Lambda $. For any $\phi (B)=g\in {\cal Q}(\Lambda )$ the equality
$B^2=B ~ {\hat \cdot } ~ B =T$ is valid (see Theorem 23). Therefore,
for each unit vector $z\in {\cal D}(T)$, one gets $z\in {\cal D}(B)$
and $<Tz;z> = <Bz;Bz>\ge 0$, hence $T\ge 0$ whenever $sp (T)\ge 0$.

\par {\bf 36. Note.} The set of positive elements in ${\cal
N}({\sf A})$ forms a positive cone. Therefore, ${\cal Q}({\sf A})$
is a partially ordered real vector space relative to the partial
ordering induced by this cone. But the unit operator $I$ is not an
order unit for ${\cal Q}({\sf A})$ in the considered case.

\section{Normal operators and homomorphisms}

\par {\bf 37. Definition.} A subset $P$ of a topological space $W$
is called nowhere dense in $W$ if its closure has empty interior. A
subset $B$ in $W$ is called meager or of the first category if it is
a countable union $B=\bigcup_{j=1}^{\infty } P_j$ of nowhere dense
subsets $P_j$ in $W$.
\par There is said that the mapping ${\sf B}(\Phi ,{\bf
R})\ni f\mapsto f(T)\in L_q(X)$ with the monotone sequential
convergence property is $\sigma $-normal.

\par {\bf 38. Lemma.} {\it Suppose that $\Lambda $ is an extremely
disconnected compact Hausdorff topological space. Then each Borel
subset of $\Lambda $ differs from a unique clopen subset by a meager
set. Each bounded Borel function $g$ from $\Lambda $ into the
Cayley-Dickson algebra ${\cal A}_v$ with $card (v)\le \aleph _0$
differs from a unique continuous function $f$ on a meager set. There
exists a conjugation-preserving $\sigma $-normal homomorphism
$\theta $ from the algebra ${\cal B}(\Lambda ,{\cal A}_v)$ of
bounded Borel functions onto $C(\Lambda ,{\cal A}_v)$. Its kernel
$ker (\theta )$ consists of all bounded Borel functions vanishing
outside a meager set.}
\par {\bf Proof.} As $card (v)\le \aleph _0$ the Cayley-Dickson algebra
is separable and of countable topological weight as the topological
space relative to its norm topology, since $card (\bigcup_{n\in \bf
N} \aleph _0^n) = \aleph _0$. \par Consider the family $\cal F$ of
all subsets contained in $\Lambda $ which differ from a clopen
subset by a meager set. Take an arbitrary $V\in \cal F$ and a clopen
subset $Q$ so that $(V\setminus Q)\cup (Q\setminus V)=: P$ is
meager, i.e. $\Lambda \setminus V$ and $\Lambda \setminus Q$ differ
by this same meager set. Since $\Lambda \setminus Q$ is clopen, one
gets $(\Lambda \setminus V)\in \cal F$. In addition each open set
$U$ in $\Lambda $ belongs to $\cal F$, since $cl (U)$ is clopen and
$cl (U)\setminus U$ is nowhere dense in $\Lambda $. For a sequence $
\{ V_j, Q_j, P_j: ~ j= 0, 1, 2,...  \} $ of such sets the inclusion
is valid:
$$ [(\bigcup_{j=1}^{\infty } V_j)\setminus (\bigcup_{j=1}^{\infty }
Q_j) ] \cup [(\bigcup_{j=1}^{\infty } Q_j)\setminus
(\bigcup_{j=1}^{\infty } V_j) ] \subseteq \bigcup_{j=1}^{\infty }
P_j .$$ The set $\bigcup_{j=1}^{\infty } P_j$ is meager and the set
$\bigcup_{j=1}^{\infty } Q_j$ is open, hence $(\bigcup_{j=1}^{\infty
} V_j)\in \cal F$. Thus the family $\cal F$ contains the $\sigma
$-algebra generated by open subsets, consequently, ${\cal F}$
contains the Borel $\sigma $-algebra ${\cal B}(\Lambda )$ of all
Borel subsets in $\Lambda $. \par Theorem 3.9.3 \cite{eng} tells
that the union $\bigcup_{j=1}^{\infty } K_j$ of a sequence of
nowhere dense subsets $K_j$ in a $\check{C}$hech-complete
topological space $W$ is a co-dense subset, i.e. its complement set
$W\setminus \bigcup_{j=1}^{\infty } K_j$ is everywhere dense in $W$.
On the other hand, each topological space metrizable by a complete
metric is $\check{C}$hech-complete, also each locally compact
Hausdorff space is $\check{C}$hech-complete \cite{eng}.
\par This implies that the complement of a meager set is dense in
$\Lambda $. Therefore, two continuous functions agree on the
complement of a meager set only if they are equal. This means that
there exists at most one continuous function agreeing with a given
bounded Borel function on the complement of a meager set.
\par Let now $V$ be a Borel subset of $\Lambda $, let also $g=\chi
_V$ be the characteristic function of $V$. Take a clopen subset $Q$
in $\Lambda $ so that $(Q \setminus V)\cup (V\setminus Q) =: P$ is a
meager subset in $\Lambda $. The characteristic function $f=\chi _Q$
of $Q$ is continuous and $(g-f)$ is zero on $\Lambda \setminus P$.
Therefore, there exists at most one clopen subset in $\Lambda $
differing from $V$ by a meager set. Therefore, a step function being
a finite ${\cal A}_v$ vector combination of characteristic functions
of disjoint Borel subsets in $\Lambda $ differs from a continuous
function on the complement of a meager set. \par The set of step
functions is dense in the ${\cal A}_v$ vector space ${\cal
B}(\Lambda ,{\cal A}_v)$ relative to the supremum-norm. This means
that if $g$ is a Borel function $g: \Lambda \to {\cal A}_v$, $~ \| g
\| := \sup_{x\in \Lambda } |g(x)|<\infty $, then there exists a
sequence of Borel step functions $\mbox{}_ng$ so that $\lim_{n\to
\infty } \| \mbox{}_ng - g \| =0$. Let $\mbox{}_nf$ be a sequence of
continuous functions agreeing with $\mbox{}_ng$ on the complement of
a meager set $P_n$. Therefore, the inequality $\| \mbox{}_nf -
\mbox{ }_mf \| \le \| \mbox{}_ng - \mbox{ }_mg \| $ if fulfilled,
since $\mbox{}_nf - \mbox{ }_mf$ and $\mbox{}_ng - \mbox{ }_mg$
agree on the complement of a meager set $P_n\cup P_m$. The set
$P_n\cup P_m$ is meager and $| \mbox{}_nf(x) - \mbox{ }_mf (x) | \le
\| \mbox{}_ng - \mbox{ }_mg \| $ for each $x\in \Lambda \setminus
(P_n\cup P_m)$. Therefore, $ \{ \mbox{}_nf : n \} $ is a Cauchy
sequence converging in supremum norm to a continuous function $f\in
C(\Lambda ,{\cal A}_v)$. Therefore, $f$ and $g$ agree on $\Lambda
\setminus (\bigcup_{n=1}^{\infty } P_n)$, where a set
$(\bigcup_{n=1}^{\infty } P_n)$ is meager. \par If functions $g^1$
and $g^2\in {\cal B}(\Lambda ,{\cal A}_v)$ differ from $f^1$ and
$f^2\in C(\Lambda ,{\cal A}_v)$ on meager sets $P^1$ and $P^2$, then
${\tilde g}^1$, $~ bg^1+g^2$, $~g^1b+g^2$ and $g^1g^2$ differ from
${\tilde f}^1$, $~bf^1+f^2$, $~f^1b+f^2$ and $f^1f^2$ on a subset of
$P^1\cup P^2$. Thus the mapping $\theta : {\cal B}(\Lambda ,{\cal
A}_v)\ni g\mapsto f\in C(\Lambda ,{\cal A}_v)$ is a conjugate
preserving surjective homomorphism so that $\theta (g)=0$ if and
only if $g$ vanishes on the complement of a meager set. \par
Particularly, if $ \{ \mbox{}_ng: ~ n \} $ is a monotone increasing
sequence in ${\cal B}(\Lambda ,{\bf R})\hookrightarrow {\cal
B}(\Lambda ,{\cal A}_v)$ of bounded Borel real-valued functions
tending point-wise to the bounded Borel function $g$, each
continuous real-valued function $\mbox{}_nf \in C(\Lambda ,{\cal
A}_v)$ differs from $\mbox{}_ng$ on the meager set $P_n$, then
$\mbox{}_nf(x)\le \mbox{ }_{n+1}f(x)$ for each $x\in \Lambda
\setminus (P_n\cup P_{n+1})$ in a dense set so that $\mbox{}_nf \le
\mbox{ }_{n+1}f$. Therefore, $\mbox{}_nf\le f$ for each natural
number $n$, while $f$ differs from $g$ on the meager set $P$. Thus
the sequences $\{ \mbox{}_ng(x): n \} $ and $\{ \mbox{}_nf(x): n \}
$ tend to $f(x)$ for each $x\in \Lambda \setminus (P\cup
(\bigcup_{n=1}^{\infty }P_n))$. This means that the function $f$ is
the least upper bound in $C(\Lambda ,{\bf R})$ of the sequence $ \{
\mbox{}_nf: n \} $ and the homomorphism $\theta $ is $\sigma
$-normal.

\par {\bf 39. Corollary.} {\it Suppose that $U$ is an open dense subset
in an extremely disconnected compact Hausdorff topological space
$\Lambda $ and $f$ is a continuous bounded function $f: U \to {\cal
A}_v$, $~card (v)\le \aleph _0$. Then there exists a unique
continuous function $\xi : \Lambda \to {\cal A}_v$ extending $f$
from $U$ onto $\Lambda $.}
\par {\bf Proof.} Put $g(x)=f(x)$ for each $x\in U$, while $g(x)=0$
for each $x\in \Lambda \setminus U$, hence $g$ is a bounded Borel
function on $\Lambda $. From Lemma 38 we know that a unique
continuous function $\xi : \Lambda \to {\cal A}_v$ exists so that
$\xi $ and $g$ differ on the complement of a meager set. If $\xi
(x)\ne f(x)$ for some $x\in U$, then by continuity of $f-\xi $ on
$U$ we get that $f(x)\ne \xi (x)$ on a non-void clopen subset $W$ of
$U$. But $W$ is not meager. This contradicts the choice of $\xi $.
Thus $\xi $ is the continuous extension of $f$ from $U$ on $\Lambda
$.

\par {\bf 40. Remark.} If $\Lambda $ is metrizable, the condition
$card (v)\le \aleph _0$ can be dropped in Lemma 38 and Corollary 39
in accordance with \S 31 in chapter 2 volume 1 \cite{kuratb}. We
denote by ${\cal B}_u(\Lambda ,{\cal A}_v)$ the family of all Borel
functions $f: \Lambda \to {\cal A}_v$.
\par {\bf 41. Lemma.} {\it Suppose that $\Lambda $ is an extremely
disconnected compact Hausdorff space. Let either $\Lambda $ be
metrizable or $card (v)\le \aleph _0$. Then a conjugation-preserving
surjective homomorphism $\psi : {\cal B}_u(\Lambda ,{\cal A}_v) \to
{\cal N}(\Lambda ,{\cal A}_v)$ exists so that its kernel $ker (\psi
) $ consists of all those functions in ${\cal B}_u(\Lambda ,{\cal
A}_v)$ vanishing on the complement of a meager set.}
\par {\bf Proof.} As $f$ and $g$ are normal functions defined on
$\Lambda \setminus W_f$ and $\Lambda \setminus W_g$ respectively and
$f(x)=g(x)$ for each $x\in \Lambda \setminus (W_f\cup W_g\cup P)$,
where $P$ is a meager subset of $\Lambda $, then $W_f=W_g$ and $f=g$
in accordance with Lemma 9, since the set $(W_f\cup W_g\cup P)$ is
meager in $\Lambda $, while $\Lambda \setminus (W_f\cup W_g\cup P)$
is dense in $\Lambda $. This implies that at most one normal
function can be agreeing with any function on the complement of a
meager set.
\par For each Borel function $g: \Lambda \to {\cal A}_v$ and each
ball $B({\cal A}_v,y,q) := \{ z\in {\cal A}_v: ~ |z-y|\le q \} $
with center $y$ and radius $q>0$, its inverse image $g^{-1}(B({\cal
A}_v,y,q))$ is a Borel subset $V_q$ of $\Lambda $. In view of Lemma
38 a clopen subset $Q_n$ exists so that $(V_n\setminus Q_n)\cup
(Q_n\setminus V_n)$ is meager in $\Lambda $. Take the Borel function
$\mbox{}_ng$ such that $\mbox{}_ng$ is equal to $g$ on $V_n$ and is
zero on $\Lambda \setminus V_n$, consequently, it satisfies the
inequality $ \| \mbox{}_ng \| \le n$. Applying Lemma 38 one gets a
continuous function $\mbox{}_nf: \Lambda \to {\cal A}_v$ so that
that agrees with $\mbox{}_ng$ on the complement of a meager subset
$P_n$ of $\Lambda $. \par The function $\mbox{}_nf$ vanishes on the
subset $\Lambda \setminus (V_n\cup P_n)$, since $\mbox{}_ng$
vanishes on $\Lambda \setminus V_n$. Certainly the set $\Lambda
\setminus (V_n\cup P_n)$ contains the subset $(\Lambda \setminus
Q_n)\setminus (P_n\cup (V_n\setminus Q_n)\cup (Q_n\setminus V_n))$.
The set $(P_n\cup (V_n\setminus Q_n)\cup (Q_n\setminus V_n))$ is
meager and $\mbox{}_nf$ is continuous on $\Lambda $, hence this
mapping $\mbox{}_nf$ vanishes on $\Lambda \setminus V_n$. The
inclusion $V_n\subseteq V_{n+1}$ for each natural number $n$ implies
the inequality $e_n(x)\le e_{n+1}(x)$ for each point $x$ in $\Lambda
$ outside a meager set, where $e_n := \chi _{Q_n}$ is the
characteristic function of $Q_n$. But continuity gives $e_n\le
e_{n+1}$ and hence $Q_n\subseteq Q_{n+1}$  for every $n\in \bf N$.
The functions $\mbox{}_{n+1}g$ and $\mbox{}_ng$ are consistent on
$V_n$, consequently, continuous functions $\mbox{}_{n+1}f$ and
$\mbox{}_nf$ agree on $Q_n$. Then the inequality $n\le
|\mbox{}_mf(x)|$ for each $x\in Q_m\setminus Q_n$ and $n<m$ is
fulfilled, since $n\le |\mbox{}_mg(y)|$ when $y\in V_m\setminus
V_n$. The equality $\bigcup_{n=1}^{\infty } V_n=\Lambda $ leads to
the inclusion
\par $W_f=\Lambda \setminus (\bigcup_{n=1}^{\infty } Q_n)
\subseteq \bigcup_{n=1}^{\infty } (V_n\triangle Q_n)$, \\
where $A\triangle B := (A\setminus B)\cup (B\setminus A)$ for two
sets. Therefore, the set $W_f$ is closed and meager, consequently,
$W_f$ is nowhere dense in $\Lambda $. If $f(x)=\mbox{ }_nf(x)$ for
each $x\in Q_n$ and $n\in \bf N$, then the mapping $f$ is continuous
on $\Lambda \setminus W_f$. \par For a point $x\in W_f$ there exists
a natural number $n$ so that $x\in \Lambda \setminus Q_n$. When
$y\in \Lambda \setminus (W_f\cup Q_n)$ the inequality $n\le |f(y)|$
is valid. Thus the function $f$ is normal. By the construction above
two functions $f$ and $g$ agree on the complement of $W_f\cup
(\bigcup_{n=1}^{\infty } P_n)$, where the latter set is meager. This
induces a conjugation-preserving surjective homomorphism $\psi :
{\cal B}(\Lambda ,{\cal A}_v)\ni g \mapsto f \in {\cal N}(\Lambda
,{\cal A}_v)$ with kernel consisting of those Borel functions
vanishing on the complement of a meager set in $\Lambda $.

\par {\bf 42. Proposition.} {\it  Let $\sf A$ be a quasi-commutative
von Neumann algebra over either the quaternion skew field or the
octonion algebra ${\cal A}_v$ with $2\le v\le 3$, let also
$\mbox{}_nT$ be a sequence of operators in ${\cal Q}({\sf A})$ with
upper bound $\mbox{}_0T$ in ${\cal Q}({\sf A})$. Then $ \{
\mbox{}_nT: ~ n=1,2,3,... \} $ has a least upper bound $T$ in ${\cal
Q}({\sf A})$.}
\par {\bf Proof.} An algebra $\sf A$ is $*$-isomorphic with
$C(\Lambda ,{\cal A}_v)$, where $\Lambda $ is a compact extremely
disconnected Hausdorff space (see Theorem I.2.52
\cite{ludopalglamb}). Take normal functions $\mbox{}_nf$ in ${\cal
N}(\Lambda ,{\cal A}_v)$ representing $\mbox{}_nT$. Then $B {\hat +}
\mbox{ }_1T$ is the least upper bound of $ \{ \mbox{}_nT: ~ n \} $,
when $B$ is the leat upper bound of $\{ \mbox{}_nT {\hat +} -\mbox{
}_1T: ~ n \} $. This means that without loss of generality we can
consider $\mbox{}_nT\ge 0$ for each $n\in \bf N$. If a normal
function $\mbox{}_nf$ is defined on $\Lambda \setminus W_n$, then
$W_n=(W_n)_+$, hence $W_n\subseteq W_{n+1}$ for each $n$. If
$x\notin W := \bigcup_{n=0}^{\infty } W_n$, then $\mbox{}_nf(x)$ is
defined for each natural number $n$ and the sequence $\{ \mbox{}_nf:
~ n \} $ has an upper bound $\mbox{}_0f(x)$ so that $\mbox{}_0f\le
h$ for each $h\in \psi ^{-1}(\phi (\mbox{}_0T))$ (see Proposition 32
and Lemma 41). Thus $ \{ \mbox{}_nf: ~ n \} $ converges to some
$g(x)$. Put $g$ to be zero on $W$, then $g$ is a Borel function on
$\Lambda $. Applying Lemma 41 one gets a normal function $f\in {\cal
N}(\Lambda ,{\cal A}_v)$ agreeing with $g$ on $\Lambda \setminus P$,
where $P$ is a meager subset in $\Lambda $. Thus the sequence $ \{
\mbox{}_nf(x): ~ n \} $ converges to $f(x)$ on the dense subset
$\Lambda \setminus (P\cup W)$. This means that $f$ is the least
upper bound for $ \{ \mbox{}_nf: ~ n \} $. Then an element $T$ in
${\cal Q}({\sf A})$ represented by $f\in {\cal Q}(\Lambda )$ is the
least upper bound of $ \{ \mbox{}_nT: ~ n=1,2,... \} $.

\par {\bf 43. Notes.} Using Proposition 42 we can extend our
definition of a $\sigma $-normal homomorphism to ${\cal N}({\sf
A})$, $~ {\cal N}(\Lambda ,{\cal A}_v)$, $~ {\cal B}_u({\cal
A}_v,{\cal A}_v)$ and $~ {\cal B}_u(\Lambda ,{\cal A}_v)$. In view
of Lemma 41 this homomorphism from ${\cal B}_u(\Lambda ,{\cal A}_v)$
onto ${\cal N}(\Lambda ,{\cal A}_v)$ is $\sigma $-normal. \par
Consider an increasing sequence $ \{ \mbox{}_ng: ~ n \} $ of Borel
functions on $\Lambda $ tending point-wise to the Borel function
$\mbox{}_0g$ and take the $\sigma $-normal functions $\mbox{}_nf$
corresponding to $\mbox{}_ng$. Then the sequence $\{ \mbox{}_nf: ~ n
\} $ has the least upper bound $\mbox{}_0f$ in ${\cal N}(\Lambda
,{\cal A}_v)$. Indeed, the functions $\mbox{}_ng$ and $\mbox{}_nf$
agree on the complement of a meager set $P_n$. Therefore, the limit
exists $\lim_n \mbox{ }_nf(x) =\mbox{}_0f(x)$ for each point $x\in
\Lambda \setminus (\bigcup_{n=0}^{\infty }P_n)$ and the subset
$\Lambda \setminus (\bigcup_{n=0}^{\infty }P_n)$ is dense in
$\Lambda $. If a function $h$ is an upper bound for $ \{ \mbox{}_nf:
~ n=1,2,3,...\}$, then $\mbox{}_nf(x)\le h(x)$ for each natural
number $n=1,2,3,...$ and points $x\in \Lambda \setminus
(\bigcup_{n=0}^{\infty }P_n)$ in the complement of the meager set
$(\bigcup_{n=0}^{\infty }P_n)$. Thus $\mbox{}_0f(x)\le h(x)$ for all
$x$ in the complement of this meager set, hence the mapping $h{\hat
+} -\mbox{ }_0f$ has non-negative values on a dense set. Thus
$\mbox{}_0f\le h$ and $\mbox{}_0f$ is the least upper bound in
${\cal N}(\Lambda ,{\cal A}_v)$ of the sequence $ \{ \mbox{}_nf: ~
n=1,2,3,... \} $.
\par Using Lemma 41 one can define $g(T)$ for an arbitrary
Borel function $g$ on $sp (T)$ for any normal operator $T$ in a
Hilbert space either over the quaternion skew field or the octonion
algebra ${\cal A}_v$ with $2\le v\le 3$. In accordance with Theorem
27 operators $T$, $T^*$ and $I$ generate a quasi-commutative von
Neumann algebra $\sf A$ over the algebra ${\cal A}_v$ with $2\le
v\le 3$ such that $T$ is affiliated with $\sf A$. It is known from
Theorems 2.24 and 2.28 \cite{lujmsalop} that $\sf A$ is isomorphic
with $C(\Lambda ,{\cal A}_v)$ for some extremely disconnected
compact Hausdorff topological space $\Lambda $. In accordance with
Theorem 29 there exists an isomorphism $\phi $ of ${\cal N}({\sf
A})$ onto ${\cal N}(\Lambda ,{\cal A}_v)$. As $\phi (T)$ is defined
on $\Lambda \setminus W$, the function $q$ defined to be zero on $W$
and $g\circ \phi (T)$ on $\Lambda \setminus W$ is Borel, $q\in {\cal
B}_u(\Lambda ,{\cal A}_v)$. But Lemma 41 says that a function $h$
exists so that $h\in {\cal N}(\Lambda ,{\cal A}_v)$ and $h$ agrees
with $q$ on the complement of a meager set in $\Lambda $. \par We
now define $g(T)$ as $\phi ^{-1}(h)$.  \par It may happen that $sp
(T)$ is a subset of a Borel set $V$ (in ${\cal A}_r$, for example)
and $g$ is a Borel function on $V$, then $g(T)$ will denote $g|_{sp
(T)}(T)$, where $g|_B$ denotes the restriction of $g$ to a subset
$B\subset V$.
\par As usually ${\cal B}(Y,{\cal A}_v)$ denotes the algebra of all
bounded Borel functions from a topological space $Y$ into the
Cayley-Dickson algebra ${\cal A}_v$, while ${\cal B}_u(Y,{\cal
A}_v)$ denotes the algebra of all Borel functions from $Y$ into
${\cal A}_v$.
\par {\bf 44. Theorem.} {\it  Let $\sf A$ be a quasi-commutative
von Neumann algebra over either the quaternion skew field or the
octonion algebra ${\cal A}_v$, $~ 2\le v\le 3$, so that $\sf A$ is
generated by a normal operator $T$ acting on a Hilbert space $X$
over the algebra ${\cal A}_v$. Then the mapping $\theta : g\mapsto
g(T)$ of the algebra ${\cal B}_u(sp (T),{\cal A}_v)$ into ${\cal
N}({\sf A})$ is a $\sigma $-normal homomorphism such that $\theta
(1)=I$ and $\theta (id)=T$. Moreover, the mapping $V\mapsto E(V)$ of
Borel subsets in ${\cal A}_v$ into $\sf A$ is an ${\cal A}_v$ graded
projection-valued measure on a Hilbert space $X$, where $E(V)=\chi
_V(T)$, $~ \chi _V$ denotes the characteristic function of $V\in
{\cal B}({\cal A}_v)$. Suppose that $h: {\cal A}_v \to {\cal A}_v$
is a Borel bounded function, ${\cal B}({\cal A}_v,{\cal A}_v)$
denotes the algebra over ${\cal A}_v$ of all such Borel bounded
functions, then
$$(1)\quad \| h(T) \| \le \sup _{y\in {\cal A}_v} |h(y)| =: \| h \|
\mbox{  and}$$  $$(2)\quad <h(T)x;x> = \int_{{\cal A}_v}  d\mu
_x(y).h(y)$$ for each vector $x\in X$ and every Borel bounded
function $h\in {\cal B}({\cal A}_v,{\cal A}_v)$, where $\mu
_x(V).h(y) := <E(V).h(y)x;x>$. For $f\in {\cal B}_u({\cal A}_v,{\cal
A}_v)$ a vector $x$ belongs to a domain ${\cal D}(f(T))$ if and only
if $$(3)\quad \int_{{\cal A}_v}  d\mu _x(y).|f(y)|^2 =: \| f(T)x \|
^2< \infty $$ and Formula $(2)$ is valid with $f$ in place of $h$.
\par Let $\phi $ be an extension on ${\cal N}({\sf A})$ of the
isomorphism $\psi $ of $\sf A$ with $C(\Lambda ,{\cal A}_v)$, let
also $\nu _x$ be a regular Borel measure on $\Lambda $ so that
$$(4)\quad <Bx;x> = \int_{\Lambda }  d \nu _x(t).(\phi (B))(t)$$ for
each $B\in \sf A$, then $x\in {\cal D}(Q)$ with $Q\in {\cal N}({\sf
A})$ if and only if
$$(5)\quad  \int_{\Lambda }  d \nu _x(t).|(\phi (Q))(t)|^2 =: \| Qx \| ^2
<\infty .$$ If additionally $T$ is a self-adjoint operator, its
spectral resolution is $ \{ \mbox{}_bE: ~ b\in {\bf R} \} $, where
$\mbox{}_bE = I - E((b,\infty ))$, and $x\in {\cal D}(f(T))$ if and
only if $$(6)\quad \int_{-\infty }^{\infty }  d
<\mbox{}_bE.|f(b)|^2x;x> <\infty $$ and for such vector $x\in {\cal
D}(f(T))$ in a domain of $f(T)$ the equality
$$(7)\quad <f(T)x;x> = \int _{-\infty }^{\infty }  d
<\mbox{}_bE.f(b)x ;x> $$ is valid.}
\par {\bf Proof.} In accordance with Remark 43 there is a $\sigma
$-normal homomorphism $\omega $ from ${\cal B}_u (sp (T))$ into
${\cal B}_u(\Lambda ,{\cal A}_v)$. On the other hand, the mapping
assigning $h\in {\cal N}(\Lambda ,{\cal A}_v)$ to $\omega (g)$ is a
$\sigma $-normal homomorphism, where $g\in {\cal B}_u (sp (T))$.
Therefore, the mapping $g\mapsto g(T)$ is a $\sigma $-normal
homomorphism from ${\cal B}_u({\cal A}_v,{\cal A}_v)$ into ${\cal
N}({\sf A})$. When $g=1$ is the constant unit function, the function
$h$ is also constant unit on $\Lambda $, consequently, $g(T)=I$.
When $g=id$, then $h=\phi (T)$, so that $id (T)=T$, where $id(x)=x$
for each $x\in {\cal A}_v$.
\par As $V$ is a Borel subset of ${\cal A}_v$ and $g=\chi _V$ is its
characteristic function, the identities are satisfied: $g(T)^* =
{\tilde g}(T) = g(T)$ and $g(T) = g^2(T) = g(T)^2$. This means that
$g(T)$ is an ${\cal A}_v$ graded projection $E(V)$ in $\sf A$.
Particularly, $\chi _{ \emptyset }=0$ and $E(\emptyset )=0$, $\chi
_{{\cal A}_v}\equiv 1$, $\chi _{{\cal A}_v}(T)=I$ and $E({\cal
A}_r)=I$.
\par Recall that ${\cal A}_v$ graded projection valued measures were
defined in \S \S I.2.73 and I.2.58. We use in this section the
simplified notation $E$ instead of ${\hat {\bf E}}$.
\par For any countable family $\{ V_j: ~ j \} $ of disjoint Borel subsets in ${\cal
A}_ v$ and their characteristic functions $\mbox{}_jg := \chi
_{V_j}$, their sums $\mbox{}_nh := \mbox{}_1g+...+\mbox{}_ng$ form
an increasing sequence tending point-wise to the characteristic
function $h := \chi _V$, where $V=\bigcup_{j=1}^{\infty } V_j$. Then
$ \{ \sum_{j=1}^n E(V_j): ~ n \} $ has the least upper bound $E(V)$
so that $E(V) = \sum_{j=1}^{\infty } E(V_j)$, since $g\mapsto g(T)$
is a $\sigma $-normal homomorphism. Using our notation we get:
$$<h(T)x;x> = <E(V)x;x> = \mu _x(V) = \int_{{\cal A}_v}  d \mu
_x(t).h(t)$$ hence Equation $(2)$ is valid for Borel step functions
$h$.
\par For any bounded Borel function $h$ on ${\cal A}_v$ with values
in ${\cal A}_v$ we have the inequality $$ \| \omega (h) \| \le \| h
\| = \sup_{x\in {\cal A}_r}  |h(x)|$$ and the function $f\in {\cal
N}(\Lambda ,{\cal A}_v)$ corresponding to $\omega (h)$ belongs to
$C(\Lambda ,{\cal A}_v)$ by Lemma 38. We put $f=\theta (h)$. Then we
infer, that
$$ \| h(T) \| = \| \phi ^{-1} (f) \| = \| f \| \le \| h \| ,$$ since
$ \| f \| \le \| \omega (h) \| \le \| h \| $. Each bounded Borel
function in ${\cal B}({\cal A}_v,{\cal A}_v)$ is a norm limit of
Borel step functions, consequently, Equality $(2)$ is valid for each
$h\in {\cal B}({\cal A}_v,{\cal A}_v)$.
\par We consider now a self-adjoint operator $T$ and the
characteristic function $g=\chi _{(b,\infty )}$ of $(b,\infty )$.
Then $\theta (g)$ is the characteristic function of $\phi
(T)^{-1}((b,\infty ))$ which is an open subset contained in $\Lambda
\setminus W$ and hence in $\Lambda $. Moreover, the function in
$C(\Lambda ,{\cal A}_v)$ corresponding to $\theta (g)$ is the
characteristic function $1- e_b$ of $cl [\phi (T)^{-1}((b,\infty
))]$. This means that $E((b,\infty )) = g(T)\in \sf A$ corresponds
to $1-e_b$ and from Theorem 16 one gets $\mbox{}_bE = I- E((b,\infty
))$.
\par This implies that $<(\mbox{}_cE-\mbox{}_bE)x;x> = \mu
_x((b,c])$ for any $b\le c$ and
$$\int_{-\infty }^{\infty } d<\mbox{}_bE.|f(b)|^2 x;x> = \int
 d \mu _x(b).|f(b)|^2$$ for each Borel function $f\in {\cal
B}_u({\cal A}_v,{\cal A}_v)$, consequently, the last assertions of
this theorem reduce to Formulas $(2,3)$. \par Denote by $q_n := \chi
_{V_n}$ the characteristic function of the subset $V_n :=
|f|^{-1}([0,n])$, where $f\in {\cal B}_u({\cal A}_r,{\cal A}_v)$ is
a Borel function, we put $\mbox{}_nF := q_n (T)$. Therefore,
$\mbox{}_nf(T)=f(T)\mbox{ }_nF$ due to the first part of the proof,
where $\mbox{}_nf = fq_n$, since the functions $f$ and $q_n$
commute, $fq_n=q_n ~f$. This implies that for a vector $x\in X$ we
get the equalities:
$$(6)\quad \| f(T) \mbox{ }_nF x \| ^2 = < |f_n|^2(T)x;x> = \int_{V_n}
d\mu _x(t).|f(t)|^2 \mbox{  and}$$
$$(7)\quad \| f(t)\mbox{ }_n Fx - f(t)\mbox{ }_m F x \| ^2 =
\int_{V_n\setminus V_m}  d \mu _x(t).|f(t)|^2 .$$ The sequence $ \{
\mbox{}_nF: ~ n \} $ is increasing with least upper bound $I$, since
$ \{ q_n : n \} $ is the increasing net of non-negative Borel
functions tending point-wise to $1$. Then we infer that
$\mbox{}_nFf(T)x = f(T) \mbox{ }_nF x$ for each vector $x$ in a
domain ${\cal D}(f(T))$, since $\mbox{}_nF f(T) \subseteq f(T)
\mbox{ }_nF$. This implies that the limit exists
$$\lim_n f(T) \mbox{ }_nF x = f(T)x$$ and Formula $(3)$ follows
from $(6)$.
\par Vise versa, if the integral $\int_{{\cal A}_v} d\mu
_x(t).|f(t)|^2 $ converges, then one gets a Cauchy sequence $ \{
f(T) \mbox{ }_nF x: ~ n \} $ due to Formula $(7)$ converging to some
vector in the Hilbert space $X$ over the Cayley-Dickson algebra
${\cal A}_v$. But $x\in {\cal D}(f(T))$, since $\lim_n \mbox{
}_nFx=x$ and the function $f(T)$ of the operator $T$ is closed. \par
Quite analogous demonstration leads to Formula $(5)$. We have a
non-negative measure $$\mu _x(V).1 = <E(V).1x;x> = <E^2(V).1x;x> =
<E(V).1x;E^*(V).1x> ~ \ge 0$$ for each Borel subset $V$. As $x$ is a
vector in the domain ${\cal D}(f(T))$, the function $f$ belongs to
the Hilbert space $L^2({\cal A}_v, \mu _x, {\cal A}_v)$, while
$L^2({\cal A}_v, \mu _x, {\cal A}_v)\subset L^1({\cal A}_v, \mu _x,
{\cal A}_v)$, since $\mu _x.1$ is a finite non-negative measure,
where $L^p({\cal A}_v, \mu _x, {\cal A}_v)$ with $1\le p<\infty $
denotes an ${\cal A}_v$ vector space which is the norm completion of
the family of all step Borel functions $u$ from ${\cal A}_v$ into
${\cal A}_v$, where the norm is prescribed by the formula:
$$ \| u \| := \sqrt[p]{\int_{{\cal A}_v} d\mu _x(t).|u(t)|^p} .$$
Finally one deduces that
$$<f(T)x;x> = \lim_n <f(T)\mbox{ }_nFx;x> =\lim_n \int_{{\cal A}_v}
q_n(t) d\mu _x(t).f(t) = \int_{{\cal A}_v} d\mu _x(t).f(t) .$$
\par {\bf 45. Remark.} The ${\cal A}_v$ graded projection $E(V)$ of
the preceding theorem will also be referred as the spectral ${\cal
A}_v$ graded projection for $T$ corresponding to the Borel subset
$V$ of ${\cal A}_v$, where $2\le v \le 3$.
\par {\bf 46. Theorem.} {\it  Suppose that $T$ is a normal operator
affiliated with a von Neumann algebra acting on a Hilbert space over
either the quaternion skew field or the octonion algebra ${\cal
A}_v$ with $2\le v \le 3$ and $\psi $ is a $\sigma $-normal
homomorphism of the algebra ${\cal B}_u({\cal A}_v,{\cal A}_v)$ of
Borel functions into ${\cal N}({\sf A})$ so that $\psi (1)=I$ and
$\psi (id) =T$, where $id : {\cal A}_v\to {\cal A}_v$ denotes the
identity mapping $id (t)=t$ on ${\cal A}_v$. Then $\psi (f) =
f^r(T)$ for each $f\in {\cal B}_u({\cal A}_v,{\cal A}_v)$, where
$f^r=f|_{sp (T)}$ denotes the restriction of $f$ to $sp (T)$.}
\par {\bf Proof.} This homomorphism $\psi $ is adjoint preserving,
since $\psi $ is $\sigma $-normal. Positive elements of ${\cal
B}_u({\cal A}_v,{\cal A}_v)$ have positive roots, hence $\psi $ is
order preserving. Moreover, $\psi : {\cal B}({\cal A}_v,{\cal
A}_v)\to \sf A$ and does not increase norm, since $\psi (1)=I$.
\par At first we consider the case when $T$ is bounded. Put
$\mbox{}_0g := \chi _{{\cal A}_v\setminus B}$, where $B=B({\cal
A}_v,0,2 \| T \| )$ is the closed ball in ${\cal A}_v$ with center
$0$ and radius $2 \| T \| $.  Then $0\le (2 \| T \| )^n \mbox{
}_0g\le |id|^n$ for each natural number $n=1,2,3,...$. Therefore,
$\psi (|id |^n) = T^n$ such that $0\le (2 \| T \| )^n \psi
(\mbox{}_0g)\le | T | ^n$, hence $ \| \psi (\mbox{}_0g) \| \le
2^{-n}$ for each natural number $n$, consequently, $\psi
(\mbox{}_0g)=0$.  \par As $\mbox{}_1g$ is the characteristic
function of the ball $B$, we get $\psi (\mbox{}_1g) = I$, such that
$\psi (\mbox{}_1gh) = \psi (h)$ for each Borel function $h\in {\cal
B}_u({\cal A}_v,{\cal A}_v)$. \par Consider the restriction
$\mbox{}_0h = h|_B$ of $h$ to $B$ and put $\psi ^0(\mbox{}_0h) =
\psi (h)$. Therefore, this mapping $\psi ^0$ is a $\sigma $-normal
homomorphism of ${\cal B}_u(B,{\cal A}_v)$ into ${\cal N}({\sf A})$
with the properties: $\psi ^0(\chi _B) =I$ and $\psi^0 (\mbox{}_0id)
= T$ and $\psi ^0: {\cal B}(B,{\cal A}_v)\to \sf A$. \par At the
same time it is known from the exposition presented above that
$C(B,{\cal A}_v)$ is a $C^*$-algebra over the algebra ${\cal A}_v$
with unit being the constant function $\chi _B$. Therefore, by
Theorem 1.3.17 \cite{ludopalglamb} one gets that $\psi ^0 (f) = \psi
^0 (f(\mbox{}_0id )) = f(T)$ for each continuous function $f\in
C(B,{\cal A}_v)$. It is known from \S I.3.17 that the characteristic
function $\mbox{}_1h$ of the open subset $B\setminus sp (T)$ of $B$
is the point-wise limit of an increasing sequence $ \{ \mbox{}_nf: ~
n \} $ of positive functions $\mbox{}_nf$ on $B$, consequently,
$\psi ^0(\mbox{}_1h)$ is the least upper bound in $\sf A$ of the
sequence $\{ \psi ^0(\mbox{}_nf): ~ n \} $, since the homomorphism
$\psi ^0$ is $\sigma $-normal. Each function $\mbox{}_nf$ is
continuous and vanishes on $sp (T)$, consequently, $\psi
^0(\mbox{}_nf) = \mbox{}_nf(T)$ and $\mbox{}_nf (T)=0$ and hence
$\psi ^0(\mbox{}_1h)=0$. For the characteristic function $\chi _{sp
(T)}$ of $sp (T)\subset B$ one obtains the equalities: $\psi ^0(\chi
_{sp (T)})=I$ and $\psi ^0(\chi _{sp (T)}\mbox{ }_0h) = \psi
^0(\mbox{}_0h)$ for each $\mbox{}_0h\in {\cal B}_u(B,{\cal A}_v)$.
\par If put $\psi ^1(q|_{sp (T)}) = \psi ^0(q)$ for each $q\in {\cal
B}_u(B,{\cal A}_v)$, then $\psi ^1 : {\cal B}_u(sp (T),{\cal
A}_v)\to {\cal N}({\sf A})$ is a $\sigma $-normal homomorphism so
that $\psi ^1(1)=I$ and $\psi ^1(id)=T$ and $\psi ^1({\cal B}(sp
(T),{\cal A}_v))\subset \sf A$. In view of Theorem I.3.21
\cite{ludopalglamb} applied to to the restriction $\psi ^1|_{{\cal
B}(sp (T),{\cal A}_v)}$ the equality $\psi ^1(f)=f(T)$ is valid for
each Borel bounded function $f\in {\cal B}(sp (T),{\cal A}_v)$. \par
On the other hand, each positive function $g$ is the point-wise
limit of an increasing sequence of positive functions in ${\cal
B}(sp (T),{\bf R})\subset {\cal B}(sp (T),{\cal A}_v)$. Therefore,
$\psi ^1(g)=g(T)$ for each positive Borel function $g\in {\cal
B}_u(sp (T),{\bf R})$, since the homomorphism $\psi ^1: {\cal
B}_u(sp (T),{\cal A}_v)\to {\cal N}({\sf A})$ is $\sigma $-normal.
Using the decomposition $h=\sum_j h_ji_j$ of each Borel function
$h\in {\cal B}_u(sp (T),{\cal A}_v)$ with real-valued Borel
functions $h_j$ and $h_j=h_j^+ - h_j^-$ with non-negative Borel
functions $h_j^+$ and $h_j^-$ we infer that $\psi ^1 (h)=h(T)$. This
implies that if $q\in {\cal B}_u({\cal A}_v,{\cal A}_v)$ and
$\mbox{}_0q = q|_B$ and $q^r =q|_{sp (T)}$, then $\psi (q) = \psi ^0
(\mbox{}_0q) = \psi ^1(q^r)=q^r(T).$
\par We take now an arbitrary normal operator $T\in {\cal N}({\sf A})$
and a bounding ${\cal A}_v$ graded projection $E$ for $T$ in $\sf
A$. The mapping $\phi $ posing $(Y{\hat \cdot } E)|_{E(X)}\in {\cal
N}({\sf A}E)$ acting on $E(X)$ to $Y\in {\cal N}({\sf A})$ is a
$\sigma $-normal homomorphism of ${\cal N}({\sf A})$ into ${\cal
N}({\sf A}E)$. Taking the composition of $\phi $ with $\psi $ yields
a $\sigma $-normal homomorphism $\psi ^2: {\cal B}_u({\cal
A}_v,{\cal A}_v)\to {\cal N}({\sf A}E)$ mapping $1$ onto $E|_{E(X)}$
and $id $ onto $T|_{E(X)}$. But the composition of $\phi $ with the
mapping $f\mapsto f^r(T)$ of ${\cal B}_u({\cal A}_v,{\cal A}_v)$
into ${\cal N}({\sf A}E)$ is another homomorphism. The restriction
$T|_{E(X)}$ is bounded, consequently, from the first part of this
proof we infer that
$$(\psi (f){\hat \cdot } E)|_{E(X)} = \psi ^2 (T|_{E(X)}) =
(f^r(T){\hat \cdot } E)|_{E(X)} .$$ Theorem 23 states that there
exists a common bounding ${\cal A}_v$ graded sequence $\{
\mbox{}_nE: ~ n \} $ for $T$ and $\psi (f)$ and $f^r(T)$, where $f$
is a given element of ${\cal B}_u({\cal A}_v,{\cal A}_v)$. Then we
deduce that
$$(\psi (f){\hat \cdot } \mbox{ }_nE)|_{\mbox{}_nE(X)} = (\psi (f)
\mbox{ }_nE)|_{\mbox{}_nE(X)} = (f^r(T){\hat \cdot } \mbox{
}_nE)|_{\mbox{}_nE(X)} = (f^r(T)\mbox{ }_nE)|_{\mbox{}_nE(X)}$$ so
that $\psi (f)\mbox{ }_nE= f^r(T) \mbox{ }_nE$ for each $n$.
Therefore, $\psi (f) = f^r(T)$, since $\bigcup_{n=1}^{\infty }
\mbox{ }_nE(X)$ is a core for both $\psi (f)$ and $f^r(T)$.
\par {\bf 47. Note.} The procedure of \S 43 assigning $g(T)\in \sf
A$ can be applied to ${\cal N}(\Omega ,{\cal A}_v)$, where $\sf Y$
is another quasi-commutative von Neumann algebra over either the
quaternion skew field or the octonion algebra ${\cal A}_v$ with
$2\le v \le 3$ so that $T$ is affiliated with $\sf Y$ and ${\sf
Y}\cong C(\Omega ,{\cal A}_v)$, $~\Omega $ is an extremely
disconnected compact Hausdorff topological space. The operator in
${\cal N}({\sf A})$ formed in this way is $g(T)\in {\cal N}({\sf
Y})$ by Theorem 46.
\par {\bf 48. Corollary.} {\it Let $T$ be a normal operator satisfying
conditions of Theorem 46 and let $f$ and $g$ be in ${\cal B}_u({\cal
A}_v,{\cal A}_v)$, where $2\le v\le 3$. Then
$$(1)\quad (f\circ g)(T) = f(g(T)).$$}
\par {\bf Proof.} Consider the von Neumann quasi-commutative von
Neumann algebra over the algebra ${\cal A}_v$ with $2\le v \le 3$
generated by $T$, $T^*$ and $I$, then $g(T)\in {\cal N}({\sf A})$
and $f\mapsto f\circ g$ is a $\sigma $-normal homomorphism $\phi $
of ${\cal B}_u({\cal A}_v,{\cal A}_v)$ into ${\cal B}_u({\cal
A}_v,{\cal A}_v)$. Taking the composition of $\phi $ with the
$\sigma $-normal homomorphism $h\mapsto h(T)$ of ${\cal B}_u({\cal
A}_v,{\cal A}_v)$ into ${\cal N}({\sf A})$ leads to a $\sigma
$-normal homomorphism $\xi : f\mapsto (f\circ g)(T)$ of ${\cal
B}_u({\cal A}_v,{\cal A}_v)$ into ${\cal N}({\sf A})$ with $\xi
(1)=I$ and $\xi (id) = g(T)$. Then Formula $(1)$ follows from
Theorem 46.
\par {\bf 49. Proposition.} {\it Suppose that $\psi $ is a
$\sigma $-normal homomorphism of ${\cal N}({\sf A})$ into ${\cal
N}({\sf B})$ so that $\psi (I)=I$, where $\sf A$ and $\sf B$ are von
Neumann quasi-commutative algebras over either the quaternion skew
field or the octonion algebra ${\cal A}_v$, where $2\le v\le 3$.
Then $\psi (f(T)) = f(\psi (T))$ for each $T\in {\cal N}({\sf A})$
and each $f\in {\cal B}_u({\cal A}_v,{\cal A}_v)$.}
\par {\bf Proof.} Each quasi-commutative algebra $\sf A$ over
the algebra ${\cal A}_v$ with $2\le v \le 3$ has the decomposition
${\sf A}=\bigoplus_j {\sf A}_j i_j$, where ${\sf A}_j$ and ${\sf
A}_k$ are real isomorphic commutative algebras for each
$j,k=0,1,2,...,2^v-1$. On the other hand, the minimal subalgebra
$alg_{\bf R} (i_0,i_j,i_k)$ for each $1\le j\ne k$ is isomorphic
with the quaternion skew field. Therefore, knowing restrictions
$\psi |_{{\sf A}_0\oplus {\sf A}_ji_j\oplus {\sf A}_ki_k\oplus {\sf
A}_li_l}$ for each $1\le j<k$ will induce $\psi $ on $\sf A$, where
$i_l=i_ji_k$. Indeed, homomorphisms $\psi _{j,k} := \psi |_{{\sf
A}_0\oplus {\sf A}_ji_j\oplus {\sf A}_ki_k\oplus {\sf A}_li_l}$ for
different $1\le j<k$ and $1\le j'<k'$ are in bijective
correspondence: $\psi _{j',k'}\circ \theta ^{j,k}_{j',k'} = \omega
^{j,k}_{j',k'}\circ \psi _{j,k}$, where $\theta ^{j,k}_{j',k'}: {\sf
A}_0\oplus {\sf A}_ji_j\oplus {\sf A}_ki_k\oplus {\sf A}_li_l \to
{\sf A}_0\oplus {\sf A}_{j'}i_{j'}\oplus {\sf A}_{k'}i_{k'}\oplus
{\sf A}_{l'}i_{l'}$ and $\omega ^{j,k}_{j',k'}: {\sf B}_0\oplus {\sf
B}_ji_j\oplus {\sf B}_ki_k\oplus {\sf B}_li_l \to {\sf B}_0\oplus
{\sf B}_{j'}i_{j'}\oplus {\sf B}_{k'}i_{k'}\oplus {\sf
B}_{l'}i_{l'}$ denote isomorphisms of von Neumann algebras over
isomorphic quaternion skew fields.
\par There exists the mapping $f\mapsto \psi (f(T))$ of ${\cal B}_u({\cal
A}_r,{\cal A}_v)$ into ${\cal N}({\sf B})$ which is a $\sigma
$-normal homomorphism so that $1\mapsto I$ and $id \mapsto \psi
(T)$. Applying Theorem 46 we get this assertion.
\par {\bf 50. Corollary.} {\it Let $\sf A$ be a quasi-commutative
von Neumann algebra on a Hilbert space $X$ over either the
quaternion skew field or the octonion algebra ${\cal A}_v$, where
$2\le v\le 3$. Let also $E$ be an ${\cal A}_v$ graded projection in
$\sf A$ and $T\eta \sf A$ and $f\in {\cal B}_u({\cal A}_v,{\cal
A}_v)$. Then the identity $f((T{\hat \cdot }E)|_{E(X)}) = (f(T){\hat
\cdot } E)|_{E(X)}$ is fulfilled.}
\par {\bf Proof.} Consider the mapping $B\mapsto (B{\hat \cdot
}E)|_{E(X)}$ which is a $\sigma $-normal homomorphism $\psi $ of
${\cal N}({\sf A})$ onto ${\cal N}({\sf A}E|_{E(X)})$ so that $\psi
(I)=E|_{E(X)}$. But $E|_{E(X)}$ is the identity operator on $E(X)$.
Then Proposition 49 implies that $f((T{\hat \cdot }E)|_{E(X)}) =
f(\psi (T)) = \psi (f(T))= (f(T){\hat \cdot } E)|_{E(X)}$.
\par {\bf 51. Note.} Consider a situation when $T$ and $B$ are
normal operators, where $T$ and $B$ may be unbounded operators whose
spectra are contained in the domain of a Borel function $g\in {\cal
B}_u({\cal A}_v,{\cal A}_v)$ and $g$ has an inverse Borel function
$f$, where $2\le v \le 3$. Then $g(T) = g(B)$ if and only if $T=B$.
To demonstrate this mention that if $g(T)=g(B)$, then $T=(f\circ
g)(T) = f(g(T)) = f(g(B)) = (f\circ g)(B)=B$ due to Corollary 48.
\par Observe particularly, that if $T^2=B^2$ with positive operators
$T$ and $B$, then $T=B$. This means that a positive operator has a
unique positive square root.
\par Let $y$ be a non-zero vector in $X$ so that $Ty= b y$ for some
$b\in {\cal A}_v$ for a normal operator $T$, let also $f\in {\cal
B}_u({\cal A}_v,{\cal A}_v)$ be a Borel function whose domain
contains $sp (T)$. Suppose in addition that $T$ is strongly right
${\cal A}_v$ linear, that is by our definition $T(xu)=(Tx)u$ for
each $x\in X$ and $u\in {\cal A}_v$. Take an ${\cal A}_v$ graded
projection $E$ with range $ \{ x: ~ Tx=bx \} $, which is closed
since $T$ is a closed operator. Therefore, $TE=bE$ on $X$, since
$ET\subseteq TE$, and hence $$f(T)y = [(f(T){\hat \cdot }
E)|_{E(X)}] y = f((T{\hat \cdot } E)|_{E(X)})y = f(bE|_{E(X)})y =
f(bI) y$$ in accordance with Corollary 50.
\par {\bf 52. Example.} Let $X$ be a separable Hilbert space over
either the quaternion skew field or the octonion algebra ${\cal
A}_v$ with $2\le v\le 3$ and an orthonormal basis $ \{ e_n: ~ n\in
{\bf N} \} $, let also $\sf A$ be the algebra of bounded diagonal
operators $Te_n=t_ne_n$, where $t_n \in {\cal A}_v$ for each natural
number $n\in {\bf N} = \{ 1, 2, 3,... \} $. Then ${\sf A}$ is
isomorphic with $C(\Lambda ,{\cal A}_v)$, where $\Lambda = \beta
{\bf N}$ is the Stone-$\check{C}$hech compactification of the
discrete space $\bf N$ due to Theorems 3.6.1, 6.2.27 and Corollaries
3.6.4 and 6.2.29 \cite{eng}, since ${\bf N}\subset \Lambda $ and
$\Lambda $ is extremely disconnected. Take points $s_n$
corresponding to the pure states $T\mapsto <Te_n;e_n>$ of this
algebra $\sf A$, where $n\in {\bf N}$. Then the set $ \{ s_n: ~ n
\in {\bf N} \} $ is dense in $\Lambda $, since from $<Te_n;e_n>=0$
for every $n$ it follows that $T=0$ and each continuous function $f:
\Lambda \to {\cal A}_v$ vanishing on $ \{ s_n: ~ n \} $ is zero.
Consider the characteristic function $\chi _{s_n}$ of the singleton
$s_n$. The projection corresponding to ${\cal A}_ve_n$ lies in $\sf
A$ and a continuous function corresponds to this projection $\chi
_{s_n}\in C(\Lambda ,{\cal A}_v)$. Thus $s_n$ is an open subset of
$\Lambda $. Therefore, $ \{ s_n: ~ n \} $ is an open dense subset in
$\Lambda $ and its complement $Z=\Lambda \setminus \{ s_n: ~ n \} $
is a closed nowhere dense subset in $\Lambda $. \par One can define
the function $h(s_n)=t_n\in {\cal A}_v$ with $\lim_n |t_n|=\infty $,
hence this function $h$ is normal and defined on $\Lambda \setminus
Z$. \par As $t_n=n\xi _n$ so that $\xi _n\in {\cal A}_r$ with $|\xi
_n|=1$ for each $n$ we get a normal function $f$ corresponding to an
operator $Q$ affiliated with $\sf A$ and $Qe_n=n\xi _n e_n$.
Choosing $t_n=(n^{1/4}-n)\xi _n$ we obtain a normal function $g$
corresponding to an operator $B$ affiliated with $\sf A$ so that
$Be_n=(n^{1/4}-n)\xi _ne_n$. Take the vector $x= \sum_{n=1}^{\infty
} n^{-1} z_n e_n$ with $z_n\in \{ \pm 1, ~ \pm \xi _n, ~ \pm \xi
_n^* \} $ for each $n$, then $x\in X$ and $\nu _x( \{ s_n \}
)=n^{-2}$. Thus one gets $$\int_{\Lambda } d\nu _x(b).|f(b)|^2  =
\infty \mbox{ and}$$
$$\int_{\Lambda } d\nu _x(b).|(f+g)(b)|^2  = \sum_{n=1}^{\infty }
n^{-3/2} <\infty .$$ The function $f+g$ is normal when defined on
$\Lambda \setminus Z$, since $|(f+g)(s_n)|=n^{1/4}\to \infty $ and
$f{\hat +}g$ corresponds to $Q{\hat +}B$. In view of Theorem 44 this
vector $x$ does not belong to ${\cal D}(Q)$, but $x\in {\cal
D}(Q{\hat +} B)$. Thus $Q+B\ne Q{\hat +}B$.
\par Take now $f$, $Q$ and $x$ as above and put $h(s_n)=n^{-3/4}\xi
_n$. The operator $C$ corresponding to $h$ is bounded and $hf$ is
normal when defined on $\Lambda \setminus Z$. Thus the functions
$hf$ corresponds to the product $C{\hat \cdot }Q$ and
$$\int_{\Lambda } d\nu _x(b).|(hf)(b)|^2  =
\sum_{n=1}^{\infty } n^{-3/2}<\infty ,$$ hence $x\in {\cal D}(C{\hat
\cdot }Q)$. Contrary $x\notin {\cal D}(CQ)$, since $x\notin {\cal
D}(Q)$. Thus $CQ\ne C{\hat \cdot }Q$, consequently, the operator
$CQ$ is not closed. In accordance with Section 23 this product in
the reverse order $QC$ is automatically closed.
\par {\bf 53. Note.} Let $T$ and $Q$ be two positive operators
affiliated with a quasi-commutative von Neumann algebra $\sf A$,
which acts on a Hilbert space $X$ over either the quaternion skew
field or the octonion algebra ${\cal A}_v$, where $2\le v\le 3$, so
that $\sf A$ is isomorphic with $C(\Lambda ,{\cal A}_v)$. Positive
normal functions $f$ and $g$ correspond to these operators $T$ and
$Q$ such that $f$ and $g$ are defined on $\Lambda \setminus W_f$ and
$\Lambda \setminus W_g$ respectively. Therefore, their sum $f+g$ is
defined on $\Lambda \setminus (W_f\cup W_g)$ and is normal and
corresponds to $T {\hat +} Q$, hence
$$0\le \int_{\Lambda } d\nu _x(t).|f(t)|^2 \le \int_{\Lambda }
d\nu _x(t).|f(t)+g(t)|^2  <\infty $$ for each vector $x\in {\cal
D}(T{\hat +}Q)$. This implies that $x\in {\cal D}(T)$ and
symmetrically $x\in {\cal D}(Q)$, consequently, $x\in {\cal D}(T+Q)$
and $T+Q = T{\hat +}Q$.
\par If drop the condition that $T$ and $Q$ are positive, but
suppose additionally that $W_f\cap W_g=\emptyset $, then two
disjoint open subsets $U_f$ and $U_g$ in $\Lambda $ exist containing
$W_f$ and $W_g$ correspondingly. Therefore, $cl (U_f)\subset \Lambda
\setminus U_g$ and the clopen set $cl (U_f)$ contains $W_f$, since a
Hausdorff topological space $\Lambda $ is extremely disconnected.
Thus $f$ is bounded on $X\setminus cl (U_f)$ and $g$ is bounded on
$cl (U_f)$. Take  the ${\cal A}_v$ graded projection in $\sf A$
corresponding to the characteristic function $\chi _{cl (U_f)}$ of
$cl (U_f)$. Then two operators $QE$ and $T(I-E)$ belong to $L_q(X)$.
But the operator $TE$ is closed, since $T$ is closed and $E$ is
bounded, consequently, $TE=T{\hat \cdot }E$. We deduce that
$$(T{\hat +} Q) E = TE {\hat +} QE = TE + QE = (T+Q)E \mbox{  and}$$
$$(T{\hat +} Q)(I-E) = (T+Q)(I-E).$$
Take a vector $x\in {\cal D}(T{\hat +}Q)$, then $Ex$ and $(I-E)x$
are in the domain ${\cal D}(T{\hat +}Q)$, hence $Ex$ and $(I-E)x$
are in ${\cal D}(T+Q)$, consequently, $x\in {\cal D}(T+Q)$ and
inevitably we get that $T{\hat +}Q = T+Q$.

\par {\bf 54. Proposition.} {\it Suppose that $\sf A$ is a quasi-commutative
von Neumann algebra over either the quaternion skew field or the
octonion algebra ${\cal A}_v$ with $2\le v \le 3$ and $T\eta \sf A$.
Let
$$(1)\quad P(z)=\sum_{k, s, ~ n_1+...+n_k\le n} \{
a_{s,n_1}z^{n_1}...a_{s,n_k} z^{n_k} \} _{q(2k)}$$ be a polynomial
on ${\cal A}_v$ with ${\cal A}_v$ coefficients $a_{s,m}$, where $k,
s \in \bf N$, $0\le n_l\in \bf Z$ for each $l$, $~n$ is a marked
natural number, $z^0 := 1$, $~z\in {\cal A}_v$, $q(m)$ is a vector
indicating on an order of the multiplication of terms in the curled
brackets, $a_{1,n_1}...a_{1,n_k}\ne 0$ for $n_1+...+n_k=n$ and
constants $c>0$ and $R>0$ exist so that $$(2)\quad c|z|^n\le
|\sum_{k, s, ~ n_1+...+n_k= n} \{ a_{s,n_1}z^{n_1}...a_{s,n_k}
z^{n_k} \} _{q(2k)}|$$ for each $|z|>R$. Then the operator
$$\sum_{k, s, ~
n_1+...+n_k\le n} \{ a_{s,n_1}T^{n_1}...a_{s,n_k} T^{n_k} \}
_{q(2k)}$$ is closed and equal to $P(T)$.}
\par {\bf Proof.} We have that $\sf A$ is isomorphic with $C(\Lambda
,{\cal A}_v)$, where $\Lambda $ is a Hausdorff extremely
disconnected compact topological space (see Theorem I.2.52
\cite{ludopalglamb}). An operator $T$ corresponds to some normal
function $f$ defined on $\Lambda \setminus W_f$ so that a set $W_f$
is nowhere dense in $\Lambda $. Therefore, the composite function
$P(f(z))$ is defined on $\Lambda \setminus W_f$ and normal. Thus
$P(f(z))$ corresponds to the polynomial $P(T)$ of the operator $T$.
On the other hand, a vector $x$ is in ${\cal D}(P(T))$ if and only
if the integral
$$(3)\quad \int_{\Lambda } d \nu _x(t).|P(f(t))|^2 <\infty $$
converges. Consider the sets $\Lambda _m := cl \{ t: ~ |f(t)|< m \}
$. Since $$\lim_{|z|\to \infty } |P(z)|=\infty $$ and due to
Condition $(2)$ there exists a positive number $m>0$ such that
$$(4)\quad \frac{1}{2}c|f(t)|^n\le \frac{1}{2} |\sum_{k, s, ~
n_1+...+n_k= n} \{ a_{s,n_1}f(t)^{n_1}...a_{s,n_k} f(t)^{n_k} \}
_{q(2k)}| \le |P(f(t))|$$  $\forall $ $t\in \Lambda \setminus
(W_f\cup \Lambda _m)$, consequently, $$(5)\quad \int_{\Lambda
\setminus (W_f\cup \Lambda _m)} d\nu _x(t).|f(t)|^{2n}<\infty .$$ A
measure $\nu _x$ is non-negative and finite on $\Lambda $, the
function $f$ is bounded on $\Lambda _m$, consequently, $f\in
L^{2n}(\Lambda ,\nu _x,{\cal A}_v)$, hence $f\in L^k(\Lambda ,\nu
_x,{\cal A}_v)$  for each $1\le k \le 2n$. The power $f^n$ of $f$ is
defined on $\Lambda \setminus W_f$ and is normal, hence $f^n$
represents $\overline{T^n}$. In view of Theorem 44 the inclusion
$x\in {\cal D}(\overline{T^k})$ is fulfilled for each $k=1,...,n$,
particularly, $x\in {\cal D}(\overline{T})={\cal D}(T)$. \par We
have $\overline{T}=T$. Suppose that $\overline{T^k}=T^k$ for
$k=1,...,j-1$. As $y\in {\cal D}(\overline{T^j})$, then $y\in {\cal
D}(T)$. Take a bounding ${\cal A}_v$ graded sequence $\mbox{}_lE$ in
$\sf A$ for $\overline{T^{j-1}}$ and $T$. Therefore, we deduce that
$\overline{T^j}\mbox{ }_mE = (\overline{T^{j-1}} {\hat \cdot }
T)\mbox{ }_mE = \overline{T^{j-1}} \mbox{ }_mE T\mbox{ }_mE$,
consequently, $\mbox{}_mE\overline{T^j}y = \overline{T^j}\mbox{
}_mEy=\overline{T^{j-1}} \mbox{ }_mE T\mbox{
}_mEy=\overline{T^{j-1}} \mbox{ }_mE Ty$. On the other hand, the
limits exist $\lim_m \mbox{ }_mE\overline{T^j}y =\overline{T^j}y $
and $\lim_m \mbox{ }_mET^jy =T^jy $. The operator
$\overline{T^{j-1}}$ is closed, consequently, $Ty\in {\cal D}(
\overline{T^{j-1}})$ and $\overline{T^{j-1}}Ty=\overline{T^j}y$,
hence $\overline{T^j}\subseteq \overline{T^{j-1}}T$. By our
inductive assumption $\overline{T^{j-1}} = T^{j-1}$. Therefore,
$\overline{T^j}\subseteq T^j$ and hence $\overline{T^j}= T^j$. Thus
by induction we get $\overline{T^k}= T^k$ for each $k=1,...,n$,
hence $x\in {\cal D}(T^k)$ for every $k=1,...,n$, consequently,
$x\in {\cal D}(P(T))$ and inevitably $$P(T) = \sum_{k, s, ~
n_1+...+n_k\le n} \{ a_{s,n_1}T^{n_1}...a_{s,n_k} T^{n_k} \}
_{q(2k)}.$$
\par {\bf 55. Remark.} For $v\ge 4$ the Cayley-Dickson algebra
${\cal A}_v$ has divisors of zero. Therefore, we have considered
mostly the quaternion and octonion cases $2\le v \le 3$. It would be
interesting to study further spectral operators over the
Cayley-Dickson algebras ${\cal A}_v$ with $v\ge 4$, but it is
impossible to do this in one article or book if look for comparison
on the operator theory over the complex field.


\begin{thebibliography}{399}

\bibitem{baez} J.C. Baez. "The octonions". Bull. Amer.
Mathem. Soc. {\bf 39: 2} (2002), 145-205.

\bibitem{brdeso} F. Brackx, R. Delanghe, F. Sommen.
"Clifford analysis" (London: Pitman, 1982).

\bibitem{dickson} L.E. Dickson. "The collected mathematical papers".
Volumes 1-5 (Chelsea Publishing Co.: New York, 1975).

\bibitem{danschw} N. Dunford, J.C. Schwartz.
"Linear operators" (J. Wiley and Sons, Inc.: New York, 1966).

\bibitem{emch} G. Emch. "M$\grave e$chanique quantique quaternionienne et
Relativit$\grave e$ restreinte". Helv. Phys. Acta {\bf 36} (1963),
739-788.

\bibitem{eng} R. Engelking. "General topology" (Heldermann: Berlin, 1989).

\bibitem{gilmurr} J.E. Gilbert, M.A.M. Murray.
"Clifford algebras and Dirac operators in harmonic analysis". Cambr.
studies in advanced Mathem. {\bf 26} (Cambr. Univ. Press: Cambridge,
1991).

\bibitem{girard} P.R. Girard. "Quaternions, Clifford algebras and
relativistic Physics" (Birkh\"auser: Basel, 2007).

\bibitem{guesprqa} K. G\"urlebeck, W. Spr\"ossig. "Quaternionic
analysis and elliptic boundary value problem" (Birkh\"auser: Basel,
1990).

\bibitem{guetze} F. G\"ursey, C.-H. Tze. "On the role of
division, Jordan and related algebras in particle physics" (World
Scientific Publ. Co.: Singapore, 1996).

\bibitem{jungexu10} M. Junge, Q. Xu. "Representation of certain
homogeneous Hilbertian operator spaces and applications". Invent.
Mathematicae {\bf 179: 1} (2010), 75-118.

\bibitem{kadring} R.V. Kadison, J.R. Ringrose. "Fundamentals
of the theory of operator algebras" (Acad. Press: New York, 1983).

\bibitem{kansol} I.L. Kantor, A.S. Solodovnikov.
"Hypercomplex numbers" ( Springer-Verlag: Berlin, 1989).

\bibitem{killipsimon09} R. Killip, B. Simon. "Sum rules and spectral
measures of Schr\"odinger operators with $L^2$ potentials". Annals
of Mathematics {\bf 170: 2} (2009), 739-782.

\bibitem{krausryan} R.S. Krausshar, J. Ryan. "Some conformally
flat spin manifolds, Dirac operators and automorphic forms". J.
Math. Anal. Appl. {\bf 325} (2007), 359-376.

\bibitem{kravchot} V.V. Kravchenko. "On a new approach for solving
Dirac equations with some potentials and Maxwell's sytem in
inhomogeoneous media". Operator Theory {\bf 121} (2001), 278-306.

\bibitem{kuratb} K. Kuratowski. "Topology" (Mir: Moscow, 1966).

\bibitem{ludoyst} S.V. Ludkovsky, F. van Oystaeyen.
"Differentiable functions of quaternion variables". Bull. Sci. Math.
(Paris). Ser. 2. {\bf 127} (2003), 755-796.

\bibitem{ludfov} S.V. Ludkovsky. "Differentiable functions of
Cayley-Dickson numbers and line integration". J. of Mathem. Sciences
{\bf 141: 3} (2007), 1231-1298.

\bibitem{lujmsalop} S.V. Ludkovsky.  "Algebras of operators in Banach
spaces over the quaternion skew field and the octonion algebra". J.
Mathem. Sciences {\bf 144: 4} (2008), 4301-4366.

\bibitem{lufjmsrf} S.V. Ludkovsky. "Residues of functions of octonion
variables". Far East Journal of Mathematical Sciences (FJMS), {\bf
39: 1} (2010), 65-104.

\bibitem{ludancdnb} S.V. Ludkovsky. "Analysis over Cayley-Dickson
numbers and its applications" (LAP Lambert Academic Publishing:
Saarbr\"ucken, 2010).

\bibitem{luspraaca} S.V. Ludkovsky, W. Sproessig.  "Ordered
representations of normal and super-differential operators in
quaternion and octonion Hilbert spaces". Adv. Appl. Clifford Alg.
{\bf 20: 2} (2010), 321-342.

\bibitem{ludspr} S.V. Ludkovsky, W. Spr\"ossig.
"Spectral theory of super-differential operators of quaternion and
octonion variables", Adv. Appl. Clifford Alg. {\bf 21: 1} (2011),
165-191.

\bibitem{ludspr2} S.V. Ludkovsky, W. Spr\"ossig.
"Spectral representations of operators in Hilbert spaces over
quaternions and octonions", Complex Variables and Elliptic
Equations, online, DOI:10.1080/17476933.2010.538845, 24 pages
(2011).

\bibitem{ludvhpde} S.V. Ludkovsky.  "Integration of vector
hydrodynamical partial differential equations over octonions".
Complex Variables and Elliptic Equations, online,
DOI:10.1080/17476933.2011.598930, 31 pages (2011).

\bibitem{ludcmft12} S.V. Ludkovsky. "Line integration of Dirac operators
over octonions and Cayley-Dickson algebras". Computational Methods
and Function Theory, {\bf 12: 1} (2012), 279-306.

\bibitem{ludopalglamb} S.V. Ludkovsky. "Operator algebras over Cayley-Dickson numbers"
(LAP LAMBERT Academic Publishing AG $\&$ Co. KG: Saarbr\"ucken,
2011).

\bibitem{oystaey} F. van Oystaeyen. "Algebraic geometry for
associative algebras". Series "Lect. Notes in Pure and Appl.
Mathem." {\bf 232} (Marcel Dekker: New York, 2000).

\bibitem{schafb} R.D. Schafer. "An introduction to non-associative
algebras" (Academic Press: New York, 1966).

\bibitem{zelditch09} S. Zelditch. "Inverse spectral problem for
analytic domains, II: ${\bf Z}_2$-symmetric domains". Advances in
Mathematics {\bf 170: 1} (2009), 205-269.


\end{thebibliography}
\end{document}